\documentclass[a4paper,12pt]{article}

\usepackage{amssymb,amsmath}
\usepackage{xypic}
\xyoption{all}
\usepackage{graphicx}

\newtheorem{defn}{Definition}[section]
\newtheorem{proposition}[defn]{Proposition}
\newtheorem{corollary}[defn]{Corollary}
\newtheorem{rem}[defn]{Remark}
\newtheorem{exm}[defn]{Example}
\newtheorem{lemma}[defn]{Lemma}
\newtheorem{theorem}[defn]{Theorem}
\newtheorem{notat}[defn]{Notation}
\newtheorem{newpar}[defn]{}

\newtheorem{xdefn}{Definition.}
\newtheorem{xproposition}{Proposition.}
\newtheorem{xcorollary}{Corollary.}
\newtheorem{xrem}{Remark.}
\newtheorem{xexm}{Example.}
\newtheorem{xlemma}{Lemma.}
\newtheorem{xtheorem}{Theorem.}
\newtheorem{xnotat}{Notation.}
\newtheorem{xnewpar}{\it}
\newtheorem{xproof}{{\it Proof. }}
\newtheorem{xproofof}{{\it Proof}}

\newenvironment{definition}{\begin{defn}\em}{\end{defn}}
\newenvironment{remark}{\begin{rem}\em}{\end{rem}}
\newenvironment{example}{\begin{exm}\em}{\end{exm}}

\newenvironment{proof}{\begin{xproof}\em}{\end{xproof}}

\newenvironment{newparagraph*}[1]{\begin{xnewpar}\hspace*{-1.5mm}{#1}. \rm}{\end{xnewpar}}

\newenvironment{definition*}{\begin{xdefn}\em}{\end{xdefn}}
\newenvironment{remark*}{\begin{xrem}\em}{\end{xrem}}
\newenvironment{example*}{\begin{xexm}\em}{\end{xexm}}
\newenvironment{notation*}{\begin{xnotat}\em}{\end{xnotat}}
\newenvironment{proposition*}{\begin{xproposition}}{\end{xproposition}}
\newenvironment{corollary*}{\begin{xcorollary}}{\end{xcorollary}}
\newenvironment{lemma*}{\begin{xlemma}}{\end{xlemma}}
\newenvironment{theorem*}{\begin{xtheorem}}{\end{xtheorem}}

\def\qed{\hspace{0.3cm}{\rule{1ex}{2ex}}}
\newcommand\V{\bigvee}
\newcommand\ie{i.e.}
\newcommand\eg{e.g.}
\newcommand\st{\mid}
\newcommand\cf{\textrm{cf.}}
\newcommand\tensor{\operatorname{\mathfrak T}}
\newcommand\downsegment{{\downarrow}}
\newcommand\frames{\textit{Frm}}
\newcommand\suplat{\textit{SL}}
\newcommand\wellsppqu{\textit{Qu}^{e}}
\newcommand\spp{\varsigma}
\newcommand\op[1]{\langle #1 \rangle}
\newcommand\pwset[1]{\wp(\,#1)}
\newcommand\opens{\operatorname{\mathcal{O}}}
\newcommand\topology{\operatorname{\Omega}}
\newcommand\groupoid{\operatorname{\mathcal{G}}}

\newcommand\ipi{\mathcal I}
\newcommand\lcc{\operatorname{{\mathcal L}^{\vee}}}
\newcommand\ssq{\mathit{StabQu}}
\newcommand\uiq{\mathit{Qu}}
\newcommand\Lind{\operatorname{\mathfrak{Q}}}
\newcommand\LindK{\Lind_{\rm K}}
\newcommand\LindT{\Lind_{\rm T}}
\newcommand\LindKfour{\Lind_{\rm K4}}
\newcommand\LindSfour{\Lind_{\rm S4}}
\newcommand\LindSfive{\Lind_{\rm S5}}
\newcommand\RT{R_{\rm T}}
\newcommand\sppT{\tensor_{\spp}}
\newcommand\K{\tensor_{\rm K}}
\newcommand\TT{\tensor_{\rm T}}
\newcommand\TKq{\tensor_{\rm  K4}}
\newcommand\TSq{\tensor_{\rm  S4}}
\newcommand\TSc{\tensor_{\rm  S5}}
\newcommand\modlatcat{\mathit{Lat}^\lozenge}
\newcommand\bimodlatcat{\mathit{Lat}^{\lozenge\blacklozenge}}

\newcommand\Tmodlatcat{\mathit{Lat}^\lozenge_{\rm T}}
\newcommand\Tbimodlatcat{\mathit{Lat}^{\lozenge\blacklozenge}_{\rm T}}
\newcommand\Kqmodlatcat{\mathit{Lat}^\lozenge_{\rm K4}}
\newcommand\Kqbimodlatcat{\mathit{Lat}^{\lozenge\blacklozenge}_{\rm K4}}
\newcommand\Sqmodlatcat{\mathit{Lat}^\lozenge_{\rm S4}}
\newcommand\Sqbimodlatcat{\mathit{Lat}^{\lozenge\blacklozenge}_{\rm S4}}
\newcommand\Scmodlatcat{\mathit{Lat}^\lozenge_{\rm S5}}
\newcommand\Scbimodlatcat{\mathit{Lat}^{\lozenge\blacklozenge}_{\rm S5}}
\newcommand\BK{B_{\rm K}}
\newcommand\BT{B_{\rm T}}
\newcommand\BKq{B_{\rm K4}}
\newcommand\BSq{B_{\rm S4}}
\newcommand\BSc{B_{\rm S5}}
\newcommand\Idl{\operatorname{Idl}}
\newcommand\bydef{\textrm{def}}
\newcommand\dom{\operatorname{dom}}
\newcommand\EX{\operatorname{EX}}
\newcommand\EF{\operatorname{EF}}
\newcommand\EG{\operatorname{EG}}
\newcommand\AX{\operatorname{AX}}
\newcommand\AF{\operatorname{AF}}
\newcommand\AG{\operatorname{AG}}
\newcommand\iter{*}
\newcommand\inv{*}

\begin{document}

\title{An algebraic generalization of Kripke structures\thanks{Research supported in part by Funda\c{c}\~ao para a Ci\^encia e a Tecnologia through program POCI 2010/FEDER and project POCI/MAT/55958/2004.}}

\author{\sc S\'{e}rgio Marcelino and Pedro Resende}

\date{}

\maketitle

\begin{abstract}
The Kripke semantics of classical propo\-sitio\-nal normal modal logic is made algebraic via an embedding of Kripke structures into the larger class of pointed stably supported quantales. This algebraic semantics subsumes the traditional algebraic semantics based on lattices with unary operators, and it suggests natural interpretations of modal logic, of possible interest in the applications, in structures that arise in geometry and analysis, such as foliated manifolds and operator algebras, via topological groupoids and inverse semigroups. We study completeness properties of the quantale based semantics for the systems K, T, K4, S4, and S5, in particular obtaining an axiomatization for S5 which does not use negation or the modal necessity operator. As additional examples we describe intuitionistic propositional modal logic, the logic of programs PDL, and the ramified temporal logic CTL.
\end{abstract}

\section{Introduction}\label{introduction}

It is well known that the set $\pwset{W\times W}$ of all the binary relations on a set $W$ has the structure of a unital involutive quantale (see \S\ref{prel}). Hence, a Kripke structure $(W,R)$ as it appears in the semantics of propositional normal modal logic \cite{HuCr,HuCr2}, where $W$ is the set of possible worlds and $R\subset W\times W$ is the accessibility relation, can be regarded as an example of a \emph{pointed} unital involutive quantale $(\pwset{W\times W},R)$. This suggests a way of generalizing the notion of Kripke structure, namely in terms of a more general pointed unital involutive quantale $(Q,\alpha)$, and the purpose of this paper is to assess the usefulness of this idea with respect to the semantics of modal logic.

Not every unital involutive quantale is suitable for this purpose, and in this paper we restrict to the notion of \emph{stably supported quantale} that has been introduced in \cite{AIM}. In order to motivate this let us consider again the quantale $\pwset{W\times W}$ of binary relations on $W$. Each relation $R\subset W\times W$ has a \emph{domain} $\dom(R)=\{x\in W\st\exists_y\ (x,y)\in R\}$ and, since
the diagonal relation
\[\Delta_W=\{(x,y)\in W\times W\st x=y\}\]
is of course isomorphic to $W$, we can equivalently replace the domain of $R$ by the set
\[\spp R=\{(x,x)\in\Delta_W\st x\in\dom(R)\}\;,\]
which we refer to as the \emph{support} of $R$. This defines an operation
\[\spp:\pwset{W\times W}\to\pwset{W\times W}\]
that preserves unions and in addition satisfies the following properties, for all $R,S\subset W\times W$:
\begin{eqnarray*}
\spp R &\subset& \Delta_W\\
\spp R & \subset& RR^\inv\\
R & \subset& \spp R\, R\\
\spp(RS)& \subset&\spp R
\end{eqnarray*}
A \emph{stably supported quantale}, or simply \emph{ssq}, is defined to be a unital involutive quantale $Q$ equipped with a sup-lattice endomorphism $\spp:Q\to Q$ which satisfies these properties; that is, for all $a,b\in Q$ we have
\begin{eqnarray*}
\spp a &\le& e\\
\spp a &\le& aa^\inv\\
a &\le& \spp a\, a\\
\spp(ab)&\le&\spp a
\end{eqnarray*}
For each $a\in Q$ the element $\spp a$ is called the \emph{support of $a$}, the operation $\spp$ itself is referred to as the \emph{support of $Q$}, and the set $\spp Q=\{\spp a\st a\in Q\}$ is necessarily a locale whose binary meet operation coincides with the multiplication of the quantale: $a\wedge b=ab$. For instance, if $Q=\pwset{W\times W}$ we have $\spp Q=\pwset{\Delta_W}\cong\pwset W$.

The formulas of propositional modal logic can be easily interpreted on any pointed ssq $(Q,\alpha)$: there should be a \emph{valuation map} $ v$ assigning to each formula $\varphi$ an element $ v(\varphi)\in\spp Q$. The properties that such a map must satisfy are clear. For instance, conjunction of formulas should be interpreted as multiplication in $Q$: $ v(\varphi\wedge\psi)= v(\varphi) v(\psi)$; and the remaining propositional connectives are equally straightforward, both classically and intuitionistically (see \S\ref{logica}).
As regards the modal operator $\lozenge$ of possibility, we impose
\[ v(\lozenge\varphi)=\spp(\alpha  v(\varphi))\;.\]
This is easily seen to yield the usual interpretation on Kripke structures: $ v(\lozenge\varphi)$ corresponds to the domain of the relation $\alpha  v(\varphi)$, where $\alpha$ is the accessibility relation.

\paragraph{Algebraic semantics.}
If $(Q,\alpha)$ is a pointed ssq then $\spp Q$ becomes a locale equipped with unary operators in a natural way (\S\S\ref{logica}, \ref{sec:lind}), and it will be seen in this paper that from any such locale it is possible to obtain a pointed ssq, in fact giving us an adjunction between two categories, and generalizing the way in which a binary relation on $W$ (\ie, an accessibility relation) corresponds classically to a unary operator $\lozenge:\pwset W\to\pwset W$.
This leads to a semantics that subsumes the classical algebraic semantics based on lattices with operators.

The systems of modal logic that are characterized by special kinds of accessibility relations can be characterized by subcategories of the category $\ssq_*$ of pointed ssqs. For instance, the system S5, which is characterized by accessibility being an equivalence relation, will correspond to the full subcategory of $\ssq_*$ whose objects $(Q,\alpha)$ are those satisfying $\alpha^2=\alpha^\inv=\alpha\ge e$. Moreover, it is easy to address multi-modal logics, such as dynamic logic, logics of time and space, etc., in terms of the category whose objects are ssqs equipped with more than one point, the various points satisfying suitable relations.

The category $\ssq$ of ssqs has pleasant properties \cite{AIM}. For instance, there is at  most one support on any unital involutive quantale, and any homomorphism of unital involutive  quantales $h:Q\to Q'$  between ssqs $Q$ and $Q'$ automatically preserves the support (hence, being stably supported is  a property rather than extra structure). The category $\ssq$ is therefore a full subcategory of the  category of  unital involutive quantales, and in fact it is reflective. There are presentations of ssqs by generators and relations, which play the role of ``Lindenbaum algebras'' in the context of the quantale-based semantics. For instance, the ``Lindenbaum quantale'' for S5 is the pointed ssq $\LindSfive$ generated by the usual Lindenbaum algebra of S5 (with defining relations ensuring that it is a bounded sublattice of $\spp\LindSfive$) with the distinguished point $\boldsymbol{\alpha}\in\LindSfive$ being subject to  the relations $\boldsymbol{\alpha}^2=\boldsymbol{\alpha}^\inv=\boldsymbol{\alpha}\ge e$. The usual Kripke structure based models for S5 can be identified with the \emph{relational representations} of $\LindSfive$, in other words the homomorphisms
\[h:\LindSfive\to\pwset{W\times W}\]
of unital involutive quantales (see  \cite{MuRe} for some properties of such representations). The same identification of models with the relational representations of a quantale applies to many other well known examples of modal logics.

\paragraph{Groupoids and inverse semigroups.}
An important aspect of the semantics described in this paper is that there are plenty of examples of ssqs besides the quantales of binary relations \cite{AIM}, arising from various geometric or analytic structures, and thus we are provided with a uniform way of defining semantic interpretations of propositional normal modal logic based on such structures. More precisely,
let $G$ be a groupoid (\ie, a small category all of whose arrows are isomorphisms). Writing $G$ also for the set of groupoid arrows, $\pwset{G}$ is an ssq; the quantales of binary relations $\pwset{W\times W}$ are precisely the quantales that arise from the so-called pair groupoid of $W$, which has $W$ as set of objects and $W\times W$ as set of arrows, with the two projections $W\times W\to W$ being the domain and codomain maps of the groupoid.
Even more generally, the topology $\topology(G)$ of any topological \'{e}tale groupoid $G$ is a sub-ssq of $\pwset G$; and for localic \'{e}tale groupoids, too, there is an ssq $\opens(G)$ associated to each groupoid $G$.

These facts are a part of the close relation \cite{AIM} between the notions of \'{e}tale groupoid (either topological or 
localic), inverse semigroup, and quantale, which can be summarized in the following diagram whose arrows denote bijections of objects up to isomorphism, or even, in the case of $\ipi$ and $\lcc$, equivalences of categories:
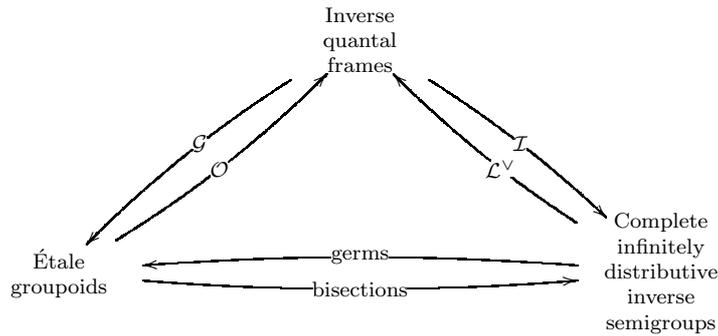
\begin{equation*}
\vcenter{\xymatrix{
&&\begin{minipage}{1.7cm}\begin{center}\scriptsize Inverse quantal frames\end{center}\end{minipage}\ar@/_/[ddll]|{\groupoid}\ar@/^/[ddrr]|{\ipi}\\
&&&&~\\
\begin{minipage}{2cm}\begin{center}\scriptsize \'{E}tale groupoids\end{center}\end{minipage}\ar@/_/[rrrr]|{\textrm{bisections}}\ar@/_/[rruu]|{\opens}&&&&
\begin{minipage}{2cm}\begin{center}\scriptsize Complete infinitely distributive inverse semigroups\end{center}\end{minipage}\ar@/^/[uull]|{\lcc}\ar@/_/[llll]|{\textrm{germs}}}} \label{trianglediagram}
\end{equation*}
The \emph{inverse quantal frames} are the ssqs that arise from \'{e}tale groupoids.

It follows that the quantale semantics automatically provides a bridge between modal logic and those areas of mathematics where examples of \'{e}tale groupoids and inverse semigroups occur, such as operator algebras and differential topology --- see, \eg, \cite{Lawson,MM,Paterson}. As an example of the latter, foliated manifolds can be associated to dynamical systems, and from a foliation it is always possible to construct a topological \'{e}tale groupoid~\cite{MM}. We shall not deal with any such examples in this paper, but we mention that if we replace $\pwset{W\times W}$ by a more general groupoid quantale, hence taking as  models of propositional modal logic the homomorphisms $\mathfrak Q\to\topology(G)$ or $\mathfrak Q\to\opens(G)$ instead of $\mathfrak Q\to\pwset{W\times W}$ (where $\mathfrak Q$ is a Lindenbaum quantale), we are led in a natural way to semantics which may be interesting, say, for applied logicians or computer scientists dealing with hybrid systems, logics of real time and space, etc.

\paragraph{Overview.}

Let $(Q,\alpha)$ be a pointed ssq. Then the locale $\spp Q$ is canonically equipped with the two unary sup-latice endomorphisms $\lozenge$ and $\blacklozenge$ defined by, for each $x\in \spp Q$,
\begin{eqnarray}
\lozenge x &=& \spp(\alpha x) \label{conj1}\\
\blacklozenge x &=& \spp(\alpha^\inv x)\label{conj2} \;,
\end{eqnarray}
which are easily seen to satisfy the following \emph{conjugacy conditions} (see \S\ref{sec:lind}):
\begin{eqnarray*}
\lozenge x\wedge y&\le&\lozenge(x\wedge\blacklozenge y)\\
\blacklozenge x\wedge y&\le&\blacklozenge(x\wedge\lozenge y)\;.
\end{eqnarray*}
Such a structure $(L,\lozenge,\blacklozenge)$, where $L$ is a locale and $\lozenge$ and $\blacklozenge$ satisfy the conjugacy conditions, will be called a \emph{bimodal frame}, and there is a functor from the category of pointed ssqs to the obvious category of bimodal frames.

A functor in the opposite direction can be easily obtained just from the knowledge  that ssqs can be presented by generators and relations; given a bimodal frame $(L,\lozenge,\blacklozenge)$ its associated quantale $Q$ is generated by the elements of $L$ plus an element $\alpha$, with relations imposing both that $L$ is a unital involutive subquantale of $Q$ and that (\ref{conj1}) and (\ref{conj2}) are satisfied. As we shall see, this defines a functor which is left adjoint to $\spp$, and in fact it is a coreflection; hence, we always have an isomorphism of bimodal frames $L\cong\spp Q$.

We can interpret this isomorphism in a logical sense as saying that no theorems are added in the process of interpreting $\lozenge$ and $\blacklozenge$ in terms of the quantale operations. Hence, if we think of the addition of the quantale operations as a language extension then this extension is conservative --- the conjugacy conditions are a complete axiomatization for the modal operators induced by $\alpha$ and $\alpha^\inv$.

This of course suggests looking at several systems of modal logic and their completeness theorems, which we shall do for K, T, K4, and S4, showing that the usual conditions on $\lozenge$ (and here also on $\blacklozenge$), as taken from the standard completeness theorems of modal logic, are precisely what is required for a coreflection to be obtained when the pointed ssqs $(Q,\alpha)$ under consideration satisfy the expected conditions, such as ``reflexivity'', ``transitivity'', etc., of the point $\alpha$. Hence, the theorems which we prove in this paper can be regarded as an algebraic generalization of the standard completeness theorems for these systems of propositional normal modal logic. It is worth noting that, as opposed to the classical theorems, these are now independent of the axiom of choice (which is required in the proof of the classical theorems, in the form of Zorn's lemma).
From our results we also obtain a complete axiomatization of system S5, where to the axioms of S4 one adds the axiom scheme
\[\lozenge \varphi\land\psi\rightarrow\lozenge(\varphi\land\lozenge \psi)\;.\]
This is obtained from the conjugacy conditions by making ${\lozenge}={\blacklozenge}$ and, as opposed to the usual axiom scheme
$\varphi\rightarrow\square\lozenge\varphi$, it does not mention negation or the modal necessity operator.

The remainder of this paper goes as follows. In \S\ref{prel} we provide some necessary background and preliminary results on quantales, presentations by generators and relations, etc. Then in \S\ref{logica} we describe in detail the quantale-based semantics of propositional normal modal logic, including the systems K, T, K4, S4, and S5. In order to illustrate the flexibility of this approach we provide additional examples, namely propositional intuitionistic logic, propositional dynamic logic, and the ramified temporal logic CTL.
Finally, after some technical results in \S\ref{sec:graded} about graded quantales (the quantale analogue of graded rings), in \S\ref{sec:lind} we address the adjunctions mentioned above. The mere existence of  the adjunctions is a consequence of the existence of presentations by generators  and relations; in other words, it can be phrased in terms of the existence of the Lindenbaum quantales. However, in order to obtain additional information about the adjunctions and,  in particular, in order to prove that they are coreflections, we shall need an actual construction of the Lindenbaum quantales. This will be conveniently formulated in terms of ``tensor algebras'' over bimodal frames, and it will take up most of \S\ref{sec:lind}.

\section{Preliminaries} \label{prel}

Here we describe some background on sup-lattices, locales and quantales. General references are
\cite{Johnstone,JT,Mu02,PaRo00,Rosenthal}.

\paragraph{Sup-lattices.}

We shall denote the category of \emph{sup-lattices} by $\suplat$. The objects are the complete lattices and the morphisms are the maps $f:L\to M$ that preserve arbitrary joins: for all $X\subset L$
\[f\left(\V X\right)=\V_{x\in X} f(x)\;.\]
We shall write $1_L$ or simply $1$ for the lattice unit (the greatest element) of a sup-lattice $L$, 
and $0_L$ or simply $0$ for the least element.

Let $(L_i)$ be a family of sup-lattices. Their cartesian product $\prod L_i$ is a sup-lattice with pointwise order and joins, and it is a product in the category $\suplat$. The products and the coproducts are isomorphic \cite{JT}, similarly to abelian groups (but, contrary to the latter, also in the case of infinitary coproducts). We denote by $\bigoplus_i L_i$ the categorical coproduct of a set-indexed family of sup-lattices $(L_i)$, and call it \emph{direct sum}.
 
The \emph{tensor product} of $L$ and $M$ is denoted by $L\otimes M$, and similarly to vector spaces, it is the image of a universal bi-morphism, where a sup-lattice \emph{bi-morphism} $f:L\times M\to N$ is a map that preserves joins in each variable separately:
\begin{eqnarray*}
f\left(\bigvee X,y\right)&=&\bigvee_{x\in X} f(x,y)\\
f\left(x,\bigvee Y\right)&=&\bigvee_{y\in Y} f(x,y)\;.
\end{eqnarray*}
Concretely, $L\otimes M$ can be identified with the set of those subsets $I\subset L\times M$ such that 
\begin{eqnarray*}
\left(x,\bigvee Y\right)\in I&\iff&\{ x\}\times Y\subset I\\
\left(\bigvee X,y\right)\in I&\iff& X\times\{ y\}\subset I
\end{eqnarray*}
for all $x\in L$, $y\in M$, $X\subset L$, and $Y\subset M$. The universal bi-morphism $L\times M\rightarrow L\otimes M$ is defined by $(x,y)\mapsto x\otimes y$, where $x\otimes y$ is the least such set that contains the pair $(x,y)$, which is called the \emph{pure tensor} generated by $(x,y)$.
$\suplat$ is a monoidal category with respect to this tensor product, with the powerset $\pwset 1$ of the singleton set as the tensor unit. Similarly to the category of abelian groups, the functor $\hom(N,-)$ has a left adjoint $-\otimes N$ for each $N$; that is, we have the familiar isomorphism 
\[\hom(M \otimes N, L) \cong \hom(M, \hom(N, L)) \, ,\] 
natural in the variables $M$ and $L$, which in fact is an order isomorphism. As a consequence of colimit preservation by left adjoints, $\otimes$ distributes over $\bigoplus$:
\[L\otimes(\bigoplus_i M_i)\cong \bigoplus_i(L\otimes M_i)\,.\]

Quotients of sup-lattices can be conveniently handled by means of closure operators (monotone endomaps $j$ that satisfy $a\le j(a)$ and $j(j(a))=j(a)$ for every element $a$). 
Let $L$ be a sup-lattice and $j$ a closure operator on $L$. The set of $j$-closed elements 
\[L_{j}=\{ x\in Q\mid x=j(x)\}\]
is a sup-lattice closed under the formation of meets in $L$, with joins given by $\bigvee^{j}(x_{i})=j(\bigvee x_{i})$,
and the map $j:L\rightarrow L_{j}$ is a (surjective) homomorphism
of sup-lattices. Conversely, given a set $S\subset L$ closed under meets in $L$, we obtain a closure operator $j_S:L\to L$ by
\[j_S(x) = \bigwedge\{y\in S\st x\le y\}\;.\]
These constructions are mutually inverse,
\[j_{L_j} = j\hspace*{3cm} L_{j_S}=S\;,\]
and every sup-lattice quotient arises in
this way up to isomorphism.

The relation to the usual description of quotients by means of congruence relations (\ie, equivalence relations on $L$ which are sub-sup-latttices of $L\times L$) is the following: from a closure operator $j$ we obtain the congruence relation $\theta_j\subset L\times L$ defined by
\[(x,y)\in\theta_j\iff j(x)=j(y)\]
[in particular $(x,j(x))\in\theta_j$] and from a congruence relation $\theta\subset L\times L$ we define a closure operator $j_\theta$ by
\[j_\theta(x) = \V[x]_\theta\]
where $[x]_\theta$ is the congruence class of $x$; of course, we have $j_{\theta_j}=j$ and $\theta_{j_\theta}=\theta$.

\paragraph{Stably supported quantales.}

A \emph{quantale} is a sup-lattice equipped with an associative \emph{multiplication}, usually written $(a,b)\mapsto ab$, which distributes over arbitrary joins:
\begin{eqnarray*}
a\left(\bigvee_{i}b_{i}\right)=\bigvee_{i}(a b_{i})\;,\hspace*{2cm}\left(\bigvee_{i}a_{i}\right) b=\bigvee_{i}(a_{i} b)\;.
\end{eqnarray*}
Hence, a quantale is a semigroup in $\suplat$.
A quantale $Q$ is \emph{unital} if the multiplication has a unit, which
we denote by $e_{Q}$, or simply $e$.

\begin{example}
A \emph{locale}, or \emph{frame}, $L$ is a sup-lattice satisfying the following distributivity property for all $x\in L$ and $Y\subset L$:
\[x\wedge\left(\V Y\right) = \V_{y\in Y} x\wedge y\;.\]
Hence, a locale is a unital quantale whose multiplication is $\wedge$ and whose unit $e$ coincides with $1$. In particular, it is a \emph{commutative} and \emph{idempotent} quantale. A quantale is a locale if and only if it is unital with $e=1$ and it is idempotent \cite{JT}.
\end{example}

An \emph{involutive quantale} $Q$ is a quantale equipped with an \emph{involution}
\[(-)^\inv:Q\rightarrow Q\;,\] 
\ie, a join preserving operation that makes $Q$ an involutive semigroup: 
\begin{eqnarray*}
(a b)^\inv&=&b^\inv a^\inv,\\
a^{\inv\inv}&=&a.
\end{eqnarray*}
Any involutive quantale satisfies $1^\inv=1$ and, if it is unital, $e^\inv=e$.

Hence, a unital involutive quantale is an involutive monoid
\[\pwset 1\stackrel e\longrightarrow Q\stackrel m\longleftarrow Q\otimes Q\]
in the monoidal category of sup-lattices, with $ab=m(a\otimes b)$.

\begin{definition}
Let $Q$ be a unital involutive quantale. A \emph{support} on $Q$ is a sup-lattice
endomorphism $\spp:Q\rightarrow Q$ satisfying, for all $a\in Q$:
\begin{eqnarray}
\spp a &\le& e\ , \label{spp:1}\\
\spp a &\le& aa^\inv\ ,\label{spp:2}\\
a &\le& \spp a a\ .\label{spp:3}
\end{eqnarray}
\end{definition}

A \emph{supported quantale} is a unital involutive quantale equipped with a specified support. On a supported quantale the set of supports $\spp Q$ coincides with $\downsegment e = \{x \mid x \le e\}$ and it is a locale with $ab=a\wedge b$ \cite{AIM}.

\begin{definition}\label{def:stable}
A support is \emph{stable} if it satisfies  $\spp(ab)=\spp(a\spp(b))$. A quantale equipped with a specified stable support is \emph{stably supported}, or simply an \emph{ssq}.
\end{definition}

Every homomorphism of unital involutive quantales preserves the support of an ssq, and thus the category of ssqs, $\ssq$, is defined to be the full subcategory of the category of unital involutive quantales $\uiq$ whose objects are the ssqs. Moreover, if a quantale is stably supported then it can have no other support, stable or not \cite{AIM}. Hence, being stably supported should be regarded as a property of unital involutive quantales rather than extra structure. In \cite{AIM} it has also been seen that the inclusion functor $\ssq\to\uiq$ has a left adjoint (\ie, $\ssq$ is a reflective subcategory of $\uiq$).

Any locale is an ssq with trivial involution and support:
\[x^\inv = x\hspace*{2cm}\spp x = x\;.\]

\paragraph{Nuclei and quotients.} \label{quotients}
The quotients of ssqs are described in a similar way to those of sup-lattices.
We give here an overview and refer to \cite{AIM} for further details.

\begin{definition}
A (quantic) \emph{nucleus} on an ssq $Q$ is a closure operator 
\[j:Q \rightarrow Q\] that satisfies, for all $x,y \in Q$,
\begin{eqnarray*}
j(x)  j(y)&\le& j(x y)\;,\\
j(x)^\inv&\le&j(x^\inv)\;,\\
\spp(j(x))&\le& j(\spp(x))\;.
\end{eqnarray*}
\end{definition}

We remark that the second condition is equivalent to $j(x)^\inv=j(x^\inv)$.

The set of $j$-closed elements $Q_j=\{x \in Q \mid x=j(x)\}$ is an ssq with joins $\bigsqcup(x_{i})=j(\bigvee x_{i})$, multiplication $x \cdot y = j(x  y)$, with the same involution as $Q$, and support $\delta x=j(\spp x)$. The map $j:Q \to Q_j$ is a (surjective) homomorphism of ssqs:
\begin{eqnarray*}
j\left(\V_i x_i\right) &=& \bigsqcup_i j(x_i)\\
j(xy) &=& j(x)\cdot j(y)\\
j(x^\inv) &=& j(x)^\inv\\
j(\spp x) &=& \delta x\;.
\end{eqnarray*}
Furthermore, every quotient arises in this way up to isomorphism.

The set $N(Q)$ of nuclei is a complete lattice under the pointwise order, with meets being calculate pointwise: 
$j\le k\Leftrightarrow\forall_{x\in Q}(j(x)\le k(x))$, and
$\bigwedge_{\alpha}j_{\alpha}(x)=\bigwedge_{\alpha}(j_{\alpha}(x))$. 
Furthermore, we have 
$j\le k\Leftrightarrow Q_k \subset Q_j$, and the join of nuclei corresponds to intersection of the respective sets of closed
elements: $j=\bigvee_{\alpha} j_{\alpha}$ if and only if $Q_j=\bigcap_{\alpha}Q_{j_{\alpha}}$.

\begin{definition}\label{Rbarra}
Let $Q$ be an ssq, and $R\subset Q\times Q$. The \emph{supported closure} $\overline R$ of the binary relation $R$ is the smallest relation that contains $R$ and is closed for the quantale operations, \ie:
\begin{eqnarray*}
R\subset \overline R\;;\\
(y,z) \in \overline R &\Rightarrow& (ay,az) \in \overline R, \textrm{ for all }  a \in Q\;;\\
(y,z) \in \overline R &\Rightarrow& (ya,za) \in \overline R, \textrm{ for all }  a \in Q\;;\\
(y,z) \in \overline R &\Rightarrow& (\spp y,\spp z) \in \overline R\;;\\
(y,z) \in \overline R &\Rightarrow& (y^\inv,z^\inv) \in \overline R\;.\\
\end{eqnarray*}
\end{definition}

Contrarily to what is done in \cite{tropo}, we shall interpret each pair $(y,z)\in R$ as an inequality $y\le z$, rather than an equation $y=z$. It is easy to see that there is a least quantic nucleus $j$ such that $j(y)\le j(z)$ for all $(y,z)\in R$: 
\begin{eqnarray*}
j_R=\bigwedge \{j \in N(Q)\mid j(y)\le j(z) \textrm{ for all } (y,z) \in R\}.
\end{eqnarray*}
Analogously to the quotients of involutive  quantales described  in \cite{tropo}, the quantale $Q_{j_R}$ has a very simple description:
\[Q_{j_R}= \{x \in Q \mid \forall_{(y,z)\in \overline R}\ (z \le x \Rightarrow y\le x)\}.\]

We stress that nuclei and quotients of unital involutive quantales equipped with any sup-lattice endomorphism $\spp:Q\to Q$ are handled in exactly the same way as described above for ssqs (we shall use this in \S\ref{sec:lind} when dealing with quantales that are just ``pre-supported''). The properties of $\spp$ pass to the quotients one by one: if $j$ is a nucleus on $Q$ and $\delta:Q_j\to Q_j$ is the sup-lattice endomorphism defined by $\delta(a)=j(\spp a)$, then if, say, the equation $\spp(ab)=\spp(a\spp b)$ holds in $Q$ then $\delta(ab)=\delta(a\delta b)$ holds in $Q_j$, etc.

\paragraph{Generators and relations.}

Let $G$ be a set (of ``generators"). The construction of the unital involutive quantale $\uiq{\op{G}}$ freely generated by $G$ is described in \cite{tropo}. Denoting by $F:\uiq\to\ssq$ the left adjoint to the inclusion $\ssq\to\uiq$ (\cf\ paragraph after \P\ref{def:stable}), it follows that $F(\uiq{\op{G}})$ is the free ssq generated by $G$, and we shall denote it by
$\ssq\op{G}$.

\begin{definition}
Let $G$ and $R\subset {\ssq{\op{G}}} \times {\ssq{\op{G}}}$ be sets. The ssq \emph{presented} by the generators in $G$ and the relations in $R$ is \[\ssq{\op{G \mid R}} \stackrel{\bydef}{=} \ssq{\op{G}}_{j_R}\;.\]
\end{definition}

If $x\in G$, one denotes by $\lbrack x \rbrack$ the image of the generator $x$ in the quantale being presented. This notation provides a useful way of describing the defining relations of a quantale presentation: we just write the conditions with respect to which the injection of generators is universal, as in the  following example for unital involutive quantales.

\begin{example}\label{exm:tensor}
Let $L$ be a sup-lattice. It follows from the universal properties of the tensor product and the direct sum of sup-lattices that the unital involutive quantale freely generated by $L$ with joins being preserved in the presentation,
\[\uiq \left\langle L \mid  \left[\V X\right] = \V_{x\in X} [x]\right\rangle\;,\]
is isomorphic to the \emph{tensor quantale}
\[TL=\bigoplus_{d\in I} L^{(d)}\] where $I$ is the free involutive monoid on one generator, whose words can be concretely identified with the strings of symbols $\alpha$ and $\alpha^\inv$ and whose unit we shall denote by $\varepsilon$, and
\[L^{(d)}=L^{\otimes \left|d\right|}=L\otimes\cdots\otimes L\ \ \ \ \ \textrm{($|d|$ times)}\;.\]
Note that $L^{(\varepsilon)}=\pwset 1$ is the neutral element of the tensor product. The multiplication is defined on pure tensors just by concatenation \[  (x_0 \otimes \ldots \otimes x_n) (y_1 \otimes \ldots \otimes y_m) = x_0 \otimes \ldots \otimes x_n \otimes y_1 \otimes \ldots \otimes y_m\;, \] where in the case of concatenation with elements of $\pwset 1$ we use the identification $\pwset 1\otimes L\cong L\cong L\otimes \pwset 1$, to produce the identity of the quantale $e=1_{\pwset 1}$. The involution
\[(-)^\inv:TL\to TL\]
is obtained from the isomorphisms $L^{(w)}\stackrel\cong\to L^{(w^\inv)}$ that are given by
\[x_1\otimes\cdots\otimes x_n\mapsto x_n\otimes\cdots\otimes x_1\;.\]
The injection of generators is the $\alpha$-coprojection of the coproduct
\[L=L^{(\alpha)}\to\bigoplus_{d\in I} L^{(d)}.\]
\end{example}

\section{Quantale semantics of modal logic} \label{logica}

\paragraph{Propositional normal modal logic.}

In this section we describe the interpretations of the classical systems of modal logic K, T, K4, S4, and S5. For details on these we refer the reader to \cite{HuCr,HuCr2}.

We shall consider fixed a set $\Pi$ of \emph{propositional symbols}. The set $\Phi$ of \emph{propositional formulas} is defined to be the least set containing $\Pi$ such that for all $\varphi,\psi\in\Phi$ we have
\begin{eqnarray*}
\neg\varphi&\in&\Phi\;,\\
\varphi\vee\psi&\in&\Phi\;,\\
\lozenge\varphi&\in&\Phi\;,
\end{eqnarray*}
where as usual we may define other connectives, for conjunction $\wedge$, implication $\rightarrow$, and the necessity modal operator $\square$, as abbreviations:
\begin{eqnarray*}
\varphi\wedge\psi &=& \neg(\neg\varphi\vee\neg\psi)\;,\\
\varphi\rightarrow\psi&=& \neg\varphi\vee\psi\;,\\
\square\varphi &=&\neg\lozenge\neg\varphi\;.
\end{eqnarray*}

\begin{definition}
A \emph{generalized Kripke model} consists of a triple $(Q,\alpha, v)$, where $Q$ is an ssq, $\alpha\in Q$ is an \emph{accessibility element}, and $ v:\Phi\to\spp Q$ is an \emph{interpretation map} satisfying the following properties for all $\varphi,\psi\in\Phi$:
\begin{eqnarray*}
 v(\varphi\vee\psi)&=&  v(\varphi)\vee v(\psi)\\
 v(\varphi) v(\neg\varphi)&=&0\\
 v(\varphi)\vee v(\neg\varphi)&=& e\\
 v(\lozenge\varphi)&=&\spp(\alpha v(\varphi))\;.
\end{eqnarray*}
\end{definition}

\begin{remark}
The above definition makes each element $v(\varphi)$ be complemented in $\spp Q$, with $v(\neg\varphi)$ being its (unique) complement, and it also follows that conjunction is interpreted as multiplication (equivalently, meet) in $\spp Q$:
\[v(\varphi\wedge\psi)= v(\varphi) v(\psi)\;.\]
This means that we interpret the formulas inside a  Boolean subalgebra of $\spp Q$, hence obtaining a classical semantics of propositional modal logic, a fact that was already implicit in the definition of the conjunction and the implication as derived connectives. However, we point out that it is easy to define a (rather natural) semantics for intuitionistic modal logic. We shall describe this at the end of \S\ref{logica}.
\end{remark}

As usual we say that a \emph{pointed} ssq consists of an ssq $Q$ together with a specified ``point" $\alpha \in Q$. A homomorphism of pointed ssqs is a homomorphism of ssqs that preserves the point:
\begin{eqnarray*}
h:(Q,\alpha)&\to&(R,\beta)\\
\alpha &\mapsto& \beta.
\end{eqnarray*}

From now on we shall denote by $\BK$ the Lindenbaum algebra of system K (\ie, the set of formulas modulo equivalence, which is a  Boolean algebra equipped with a finite join preserving endomorphism $\lozenge$).

\begin{definition}
The \emph{Lindenbaum quantale} for K is the pointed ssq $\LindK$ which is presented by generators and relations with $\BK$ as the set of generators 
and with the following relations for all $x,y\in \BK$, where we denote the selected point by $\boldsymbol{\alpha}$:
\begin{eqnarray*}
\lbrack x\vee y\rbrack &=& \lbrack x\rbrack\vee \lbrack y\rbrack\\
\lbrack \neg x\rbrack\lbrack x\rbrack &=&0\\
\lbrack \neg x\rbrack\vee\lbrack x\rbrack&=& e\\
\lbrack\lozenge x\rbrack&=&\spp(\boldsymbol{\alpha}\lbrack x\rbrack)\:.
\end{eqnarray*}
\end{definition}

From the universal property of ssqs presented by generators and relations we immediately obtain:

\begin{theorem}
There is a bijective correspondence between abstract Kripke models $(Q,\alpha,v)$ and homomorphims of unital involutive quantales
\[\LindK\longrightarrow Q\;.\]
In particular, if $W$ is a set then a homomorphism
\[\rho:\LindK\longrightarrow \pwset{W\times W}\]
is the same as a model for system K with set of possible worlds $W$ and accessibility relation $\rho(\boldsymbol{\alpha})$.
\end{theorem}

In order to obtain similar facts for other systems, such as T, K4, S4, S5, one must define the appropriate Lindenbaum quantales.

\begin{definition}\label{def:Lindquantales}
The \emph{Lindenbaum quantales} for T, K4, S4, and S5, are the pointed ssqs $\LindT$, $\LindKfour$, $\LindSfour$, and $\LindSfive$, respectively, which are presented by generators and relations similarly to $\LindK$, with the following additional relations:
\begin{description}
\item[$\LindT$:] $e\le\boldsymbol{\alpha}$
\item[$\LindKfour$:] $\boldsymbol{\alpha}\boldsymbol{\alpha}\le\boldsymbol{\alpha}$
\item[$\LindSfour$:] $e\le\boldsymbol{\alpha} \ge\boldsymbol{\alpha}\boldsymbol{\alpha}$
\item[$\LindSfive$:] $e\le\boldsymbol{\alpha}=\boldsymbol{\alpha}^\inv\ge\boldsymbol{\alpha}\boldsymbol{\alpha}$
\end{description}
\end{definition}

Hence, $\LindT$ is the quotient of $\LindK$ by the least nucleus $j$ such that $j(e)\le j(\boldsymbol{\alpha})$, 
and $\LindKfour$ is the quotient of $\LindK$ by the least nucleus $j$ such that $j(\boldsymbol{\alpha}\boldsymbol{\alpha})\le j(\boldsymbol{\alpha})$. 
Then we have $\LindSfour=\LindT\cap\LindKfour$ and $\LindSfive$ is the quotient of $\LindSfour$ by the least nucleus $j$ such that $j(\boldsymbol{\alpha})=j(\boldsymbol{\alpha}^\inv)$.

Notice that the relational representations of these quantales correspond to the expected classes of models:

\begin{theorem}
The relational representations $Q\to\pwset{W\times W}$ of the Lindenbaum quantales $Q$ correspond bijectively to the Kripke models whose accessibility relations are:
\begin{itemize}
\item Reflexive, for $Q=\LindT$;
\item Transitive, for $Q=\LindKfour$;
\item Preorders, for $Q=\LindSfour$;
\item Equivalence relations, for $Q=\LindSfive$.
\end{itemize}
\end{theorem}

\paragraph{Propositional ramified temporal logic.}

Now we describe a similar semantics for the ramified temporal logic known as \emph{Computational Tree Logic (CTL)}, see \cite{CTL}.
As above, $\Pi$ is a fixed set of \emph{propositional symbols}. The set $\Phi$ of \emph{CTL formulas} is defined to be the least set containing $\Pi$ such that for all $\varphi,\psi\in\Phi$ we have
\begin{eqnarray*}
\neg\varphi,
\varphi\vee\psi,
\EX\varphi,
\EF\varphi,
\EG\varphi\in\Phi\;,
\end{eqnarray*}
and we may define other modal operators as abbreviations: 
\begin{eqnarray*}
\AX \varphi &\equiv& \neg \EX \neg \varphi\;,\\
\AG \varphi &\equiv& \neg \EF \neg \varphi\;,\\
\AF \varphi &\equiv& \neg \EG \neg \varphi\;.
\end{eqnarray*}
\begin{definition}

The intuitive meaning of the various modalities is the following:
\begin{itemize}
\item $\EX\varphi$ means that there is a possible future where $\varphi$ will hold in the next time instant;
\item $\EF\varphi$ means that there is a possible  future where $\varphi$ will eventually hold;
\item $\EG\varphi$ means that there is a possible future where $\varphi$ will always hold (including now);
\item $\AX\varphi$ means that $\varphi$ will certainly hold in the next time instant;
\item $\AG\varphi$ means that $\varphi$ will always hold (including now) in all possible future paths;
\item $\AF\varphi$ means that in each possible future path $\varphi$ will eventually hold.
\end{itemize}

A \emph{generalized CTL model} consists of a triple $(Q,\alpha, v)$, where $Q$ is an ssq, $\alpha\in Q$ is an \emph{accessibility element} that satisfies \[\spp(\alpha)=e\]
(\ie, time never ends), and $ v:\Phi\to\spp Q$ is an \emph{interpretation map} satisfying the following properties for all $\varphi,\psi\in\Phi$:
\begin{eqnarray*}
 v(\varphi\vee\psi)&=&  v(\varphi)\vee v(\psi)\\
 v(\varphi) v(\neg\varphi)&=&0\\
 v(\varphi)\vee v(\neg\varphi)&=& e\\
 v(\EX\varphi)&=&\spp(\alpha v(\varphi))\\
 v(\EF\varphi)&=&\spp\left(\V_n\alpha^n v(\varphi)\right)\\
 v(\EG\varphi)&=&\V\{a\in Q\st a\leq v(\varphi)\wedge \spp(\alpha a)\}\;.
 \end{eqnarray*}
\end{definition}

\noindent It is easy to see that this interpretation conveys the intended meaning to the modal operators if we let $Q=\pwset{W\times W}$ for some set $W$. Only the last condition, for $\EG$, deserves an explanation. This says that $\EG\varphi$ may be interpreted as the largest subset $X\subset W$ such that every world $x\in X$ satisfies the following two conditions:
\begin{itemize}
\item $x$ satisfies the formula $\varphi$;
\item there is a world $y\in X$ such that $(x,y)\in\alpha$.
\end{itemize}
This guarantees that there is an infinite path (possibly with repetitions)
\[x_0, x_1, x_2, x_3, \ldots\]
satisfying $\varphi$ starting at any world $x_0$ where $\EG\varphi$ holds. Mathematically, the definition of $v(\EG\varphi)$ is clarified by the Knaster--Tarski fixed point theorem: the set of all the lowerbounds of $v(\varphi)$,
\[S=\{a\in\spp Q\st a\le v(\varphi) \}\;,\]
is a complete lattice and thus the set of pre-fixed points of the monotone operator
$f:S\to S$
defined by
\[f(a)= v(\varphi)\wedge\spp(\alpha\,a)\]
has a join, which in fact is a fixed point; hence, $v(\EG\varphi)$ is also the largest element $a\in\spp Q$ such that
$a=v(\varphi)\wedge\spp(\alpha\, a)$.

\paragraph{Propositional dynamic logic.}

In the program logic of \cite{PDL} there are modalities indexed by programs, which themselves form a set with some algebraic structure. Let $\Pi$ be a fixed set of \emph{propositional symbols} and $\Xi$ a set of atomic programs. The sets $F$, of \emph{formulas}, and $P$, of \emph{programs}, give us the \emph{PDL language} $\Phi=F\cup P$, and they are defined to be the least sets such that
\begin{eqnarray*}
\Pi&\subset& F\;,\\
\Xi&\subset& P\;,\\
\textrm{if }\varphi,\psi \in F &\textrm{ then }& \neg \varphi,\ \varphi\vee\psi \in F\;,\\
\textrm{if }p,q \in P &\textrm{ then }& p \cup q,\ p ; q,\ p^\iter  \in P\;,\\
\textrm{if }\varphi \in F \textrm{ and }p \in P &\textrm{ then }& \op{p}\varphi \in F\;,\\
\textrm{if }\varphi \in F &\textrm{ then }& \varphi ? \in P\;.
\end{eqnarray*}
Very briefly, the intuitive meaning of the program constructs is the following:

\begin{itemize}
\item $p\cup q$ is the program that behaves in a nondeterministic way either like $p$ or like $q$ (the choice is made at the beginning of the execution of the program, once and for all);
\item $p;q$ is the program whose execution is that of $p$ followed by $q$;
\item $p^\iter$ is the iteration of $p$, to be thought of as $p$ executed sequentially zero or more times (not to be confused with the notation for the quantale involution);
\item $\varphi ?$ is the program that tests $\varphi$, succeeding if $\varphi$ is found to be true, and failing otherwise.
\end{itemize}

\begin{definition}
A \emph{generalized PDL model} consists of a pair $(Q, v)$, where $Q$ is an ssq and $v:\Phi\to Q$ is an \emph{interpretation map} satisfying the following properties for all $\varphi,\psi\in F$ and $p, q \in P$:
\begin{eqnarray*}
 v(\varphi\vee\psi)&=&  v(\varphi)\vee v(\psi)\\
 v(\varphi) v(\neg\varphi)&=&0\\
 v(\varphi)\vee v(\neg\varphi)&=& e\\
  v(\op{p}\varphi)&=&\spp(v(p)v(\varphi))\\
v(p\cup q)&=& v(p)\vee v(q)\\
v(p;q)&=&v(p)v(q)\\
 v(p^\iter)&=&\V_{n\in\omega} v(p)^n\\
 v(\varphi ?)&=&v(\varphi)\;.
\end{eqnarray*}
\end{definition}
This interpretation shows that to a large extent both the formulas and the programs are treated on an equal footing. In particular, $p;q$ can be regarded as the (noncommutative) ``conjunction'' of $p$ and $q$, and $p\cup q$ as their disjunction, while a formula is just a particular kind of program ($\varphi$ is identified with $\varphi ?$).

\paragraph{Intuitionistic modal logic.}

It is easy to define a semantics for intuitionistic propositional modal logic
if we let $\wedge$ and $\rightarrow$ be independent connectives.
This is because the support $\spp Q$ of any ssq $Q$ is a locale and therefore a Heyting algebra, and thus, denoting by $\backslash$ the residuation operation of $\spp Q$,
\[b\backslash a = \V\{c\in\spp Q\st b\wedge c\le a\}\;,\]
the conditions on $v$ can be simply replaced by the following:
\begin{eqnarray*}
 v(\varphi\vee\psi)&=&  v(\varphi)\vee v(\psi)\\
 v(\varphi\wedge\psi)&=& v(\varphi) v(\psi)\\
v(\varphi\rightarrow\psi)&=& v(\varphi)\backslash v(\psi)\\
v(\neg\varphi)&=& v(\varphi)\backslash 0\\
  v(\lozenge\varphi)&=&\spp(\alpha v(\varphi))\;.
\end{eqnarray*}
This would entirely define the intuitionistic semantics if we contented ourselves with defining ${\square}={\neg\lozenge\neg}$ as before. 
However, this is a bad interpretation of $\square$, as for instance it usually does not satisfy the axiom of distributivity over meets
\[\square(\varphi\wedge\psi)\leftrightarrow\square\varphi\wedge\square\psi\;.\]
Indeed, a much better interpretation, in particular one that satisfies distributivity of $\square$ over (arbitrary) meets, is obtained if we let $\square$ be an independent connective interpreted as the right adjoint of a suitable sup-lattice endomorphism, as we now describe.

Let $W$ be a set, let $R\subset W\times W$ be a binary relation, and let $Q$ be the ssq $\pwset {W\times W}$. Let $\lozenge$ and $\blacklozenge$ be the sup-lattice endomorphisms of $\spp Q$ defined as in \S\ref{introduction}:
\begin{eqnarray*}
\lozenge X &=& \spp(RX)\\
\blacklozenge X &=& \spp(R^\inv X)\;.
\end{eqnarray*}
Equivalently, taking into account the isomorphism $\spp Q\cong\pwset W$ we may consider $\lozenge$ and $\blacklozenge$ to be endomorphisms of $\pwset W$:
\begin{eqnarray*}
\lozenge X &=& \{y\in W\st \exists_{x\in X}\ (y,x)\in R\}\\
\blacklozenge X &=& \{y\in W\st \exists_{x\in X}\ (x,y)\in R\}\;.
\end{eqnarray*}
It is straightforward to verify that the usual necessity operator
\[\square:\pwset W\to\pwset W\;,\]
which is defined by
\[\square X = \{y\in W\st\forall_{x\in W}\ (y,x)\in R\Rightarrow x\in X\}\;,\]
is right adjoint to $\blacklozenge$, and we may take this as the natural definition of $\square$ when such a ``possibility in the past''  operator $\blacklozenge$ is available --- similarly, a ``necessity in the past'' operator $\blacksquare$ can be defined to be the right adjoint of $\lozenge$:
\begin{eqnarray*}
\blacklozenge x\le y&\iff&x\le\square y\\
\lozenge x\le y&\iff&x\le\blacksquare y\;.
\end{eqnarray*}
 This leads to the following quantale-based intuitionistic semantics for propositional modal logic, where we assume that $\wedge$, $\rightarrow$, and $\square$ are independent connectives:

\begin{definition}
A \emph{generalized intuitionistic Kripke model} consists of a triple $(Q,\alpha, v)$, where $Q$ is an ssq, $\alpha\in Q$ is an \emph{accessibility element}, and $ v:\Phi\to\spp Q$ is an \emph{interpretation map} satisfying the following properties for all $\varphi,\psi\in\Phi$:
\begin{eqnarray*}
 v(\varphi\vee\psi)&=&  v(\varphi)\vee v(\psi)\\
 v(\varphi\wedge\psi)&=& v(\varphi) v(\psi)\\
v(\varphi\rightarrow\psi)&=& v(\varphi)\backslash v(\psi)\\
v(\neg\varphi)&=& v(\varphi)\backslash 0\\
  v(\lozenge\varphi)&=&\spp(\alpha v(\varphi))\\
  v(\square\varphi) &=& \V\{x\in\spp Q\st \spp(\alpha^\inv x)\le v(\varphi)\}\;.
\end{eqnarray*}
\end{definition}

This definition illustrates a canonical way in which to define intuitionistic semantics for other propositional modal logics, including all the examples seen earlier in this section. We stress the fact that involutive quantales are essential for this, since without the involution we would neither have the $\blacklozenge$ operator, nor a canonical definition of $\square$.

It is worth commenting on some aspects of the intuitionistic version of system S5, which similarly to its classical counterpart is based on imposing that the accessibility element $\alpha$ should be self-adjoint, and thus ${\blacklozenge}={\lozenge}$. The unit of the adjunction relating $\blacklozenge$ and $\square$ is the inequality
\[x\le\square\blacklozenge x\;,\]
and thus for intuitionistic S5 we conclude that the axiom-scheme
\begin{equation}\label{unitS5}
\varphi\rightarrow\square\lozenge\varphi
\end{equation}
is always satisfied. It is well known that this scheme (together with those for S4) characterizes the classical system S5. Another axiom-scheme which is always satisfied by intuitionistic S5 is
\[\lozenge\square\varphi\rightarrow\varphi\;,\]
which corresponds to the co-unit of the adjunction, and which classically (\ie, with ${\square}={\neg\lozenge\neg}$) is equivalent to (\ref{unitS5}).

\section{Graded unital involutive quantales}\label{sec:graded}

\paragraph{Basic definitions and properties.}

The usual notion of grading of a ring has a straightforward counterpart for quantales, which will be useful in \S\ref{sec:lind}. We shall study it now.

\begin{definition}
Let $M$ be an involutive monoid. A unital involutive quantale $Q$ is \emph{graded} over $M$ if there is an $M$-indexed family of sub-sup-lattices $Q^{(m)}$ of $Q$ satisfying the following two properties.
\begin{enumerate}
\item For each $a\in Q$ there is one, and only one, element \[(a_m)\in\bigoplus_{m\in M} Q^{(m)}\] such that \[a=\V_{m\in M} a_m\;.\]
\item The following conditions hold for all $m,n\in M$:
\begin{eqnarray*}
Q^{(m)} Q^{(n)} &\subset& Q^{(mn)}\\
1_{Q^{(\varepsilon)}}&=& e\\
\left(Q^{(m)}\right)^\inv &\subset& Q^{(m^\inv)}\;.
\end{eqnarray*}
(The latter is equivalent to $\left(Q^{(m)}\right)^\inv = Q^{(m^\inv)}$.)
\end{enumerate}
For each $m\in M$ the set $Q^{(m)}$ is called the \emph{component of $Q$ in degree $m$}.
\end{definition}

\begin{example}\label{exm:graded}
Recall the tensor quantale of \P\ref{exm:tensor}: if $L$ is a sup-lattice and $I$ is the free involutive monoid on one generator $\alpha$ then the tensor quantale
\[TL=\bigoplus_{d\in I} L^{(d)}\]
has an obvious grading over $I$ such that $L^{(\varepsilon)}\cong\pwset 1$
and $L^{(\alpha)}\cong L$.
\end{example}

The following properties are immediate:

\begin{proposition}
Let $Q$ be a unital involutive quantale graded over an involutive monoid $M$.
\begin{enumerate}
\item The map $(a_m)\mapsto\V_{m\in M} a_m$ is an isomorphism of sup-lattices
\[\bigoplus_{m\in M} Q^{(m)}\to Q\;.\]
\item If $m\neq n$ then $Q^{(m)}\cap Q^{(n)}=\{0\}$.
\item The union $\bigcup_{m\in M} Q^{(m)}$ is join-dense in $Q$.
\item $\downsegment Q^{(m)}=Q^{(m)}$ for all $m\in M$.
\end{enumerate}
\end{proposition}

There is a convenient alternative definition if the unital involutive quantale $Q$ is also a locale (an example is the quantale $TL$ of \P\ref{exm:graded} if $L$ is a locale, or the quantale $\tensor(L)$ of \S\ref{sec:lind}):

\begin{proposition}
Let $M$ be an involutive monoid, and let $Q$ be a unital involutive quantale which is also a locale. Then $Q$ is graded over $M$ if and only if there is an $M$-indexed family $(e^{(m)})$ of elements of $Q$ satisfying the following properties:
\begin{enumerate}
\item $1_Q=\V_{m\in M} e^{(m)}$ (\ie, $\left(e^{(m)}\right)$ covers $Q$);
\item $e^{(m)}\wedge e^{(n)}=0$ if $m\neq n$;
\item $e^{(m)}e^{(n)}\le e^{(mn)}$;
\item $e^{(\varepsilon)}=e$;
\item $\left(e^{(m)}\right)^\inv\le e^{(m^\inv)}$ (equiv., $\left(e^{(m)}\right)^\inv = e^{(m^\inv)}$).
\end{enumerate}
\end{proposition}

\begin{proof}
It is clear that if $Q$ is graded then it satisfies 1--5 if we let $e^{(m)}$ be $1_{Q^{(m)}}$ for each $m\in M$. For the converse we define the component $Q^{(m)}$ to be $\downsegment e^{(m)}$: then any element $a\in Q$ equals
\[a\wedge 1 = a\wedge \V_m e^{(m)}=\V_m a\wedge e^{(m)}\]
where $a\wedge e^{(m)}\in Q^{(m)}$ by definition of $Q^{(m)}$; and if $\V_m a_m = \V_m b_m$ then, for each $n\in M$ we have
\[a_n = a_n\wedge\V_m a_m =a_n\wedge\V_m b_m =\V_m a_n\wedge b_m=a_n\wedge b_n\;,\]
and in a similar way we obtain
$b_n = b_n\wedge a_n$. Hence, we have
$a_n=b_n$, and thus each element $a\in Q$ can be obtained uniquely as a join $\V_m a_m$. The rest is straightforward. \qed
\end{proof}

\paragraph{Graded nuclei and quotients.}

The natural notion of quotient that respects the grading of a quantale is provided by the following definition:

\begin{definition}
Let $Q$ be a unital involutive quantale graded over an involutive monoid $M$. A nucleus $j:Q\to Q$ is \emph{graded} if it satisfies the following two conditions for all $(a_m)\in\bigoplus_{m\in M} Q^{(m)}$:
\begin{enumerate}
\item $j\left(Q^{(m)}\right)\subset Q^{(m)}$;
\item $j\left(\V_{m\in M} a_m\right) = \V_{m\in M} j(a_m)$.
\end{enumerate}
\end{definition}

\begin{proposition}
Let $j$ be a graded nucleus as in the definition above. Then,
\begin{enumerate}
\item $j(0)=0$ (the nucleus is ``dense'');
\item $Q_j$ is graded, with each component being defined by \[(Q_j)^{(m)}=j\left(Q^{(m)}\right)\;.\]
\end{enumerate}
\end{proposition}

\begin{proof}
The first condition is obvious from the fact that $0\in\bigcap_{m\in M} Q^{(m)}$, and thus $j(0)\in\bigcap_{m\in M} Q^{(m)}=\{0\}$. For the second condition notice that if $a\in Q_j$ then on one hand we have a unique representation of $a$ as a join
\[a=\V_{m\in M} a_m\]
and, on the other hand,
\[a=j(a)=j\left(\V_{m\in M} a_m\right)=\V_{m\in M} j(a_m)\;,\]
and thus $a_m=j(a_m)$ for all $m\in M$; that is, the element $a_m$ is necessarily in $(Q_j)^{(m)}$. The rest is straightforward. \qed
\end{proof}

The nucleus induced by a binary relation is graded provided the relation respects the grading:

\begin{proposition}\label{prop:gradedrelation}
Let $Q$ be a unital involutive quantale graded over an involutive monoid $M$. Let also $R\subset Q\times Q$ be a binary relation that respects the grading in the sense that
$R\subset\bigcup_m Q^{(m)}\times Q^{(m)}$.
Then $j_R$ is a graded nucleus.
\end{proposition}

\begin{proof}
Let ${\sqsubset_R}\subset Q\times Q$ be the preorder defined by
\[a\sqsubset_R b \iff j_R(a)\le j_R(b)\;.\]
Since this is also a sub-involutive-quantale of $Q\times Q$, let us call it a \emph{congruence preorder}. By a simple adaptation of the comments at the end of the sup-lattices section of \S\ref{prel}, there is a bijection between congruence preorders and nuclei on $Q$, and $\sqsubset_R$ is the least congruence preorder on $Q$ which contains $R$. It is clear that $\sqsubset_R$ respects the grading because $R$ does, and thus if $a\in Q^{(m)}$ then $j_R(a)\in Q^{(m)}$, showing that $j_R$ satisfies the first of the properties of graded nuclei. In order to prove that it also satisfies the second property
let $k:Q\to Q$ be the map defined, for each $(a_m)\in \bigoplus_m Q^{(m)}$, as follows:
\[k\left(\V_m a_m\right) = \V_m j_R(a_m)\;.\]
Since $j_R$ is monotone we have $k\le j_R$:
\[j_R\left(\V_m a_m\right) \ge \V_m j_R(a_m)= k\left(\V_m a_m\right)\;.\]
Now let us see that $k$ is itself a nucleus. First, it is obvious that it is monotone and that it satisfies $a\le k(a)$ for all $a\in Q$. It is also idempotent because
\[k\left(k\left(\V_m a_m\right)\right) =
k\left(\V_m j_R(a_m)\right)\]
and the fact that $j_R(a_m)\in Q^{(m)}$ allows us to conclude that the right hand side of the above equation equals
\[\V_m j_R(j_R(a_m))=\V_m j_R(a_m)=k\left(\V_m a_m\right)\;.\]
Now let us prove the condition relating $k$ to the multiplication. For each pair $(a_m), (b_m)\in\bigoplus_m Q^{(m)}$ we have
\begin{eqnarray*}
k\left(\V_m a_m\right)k\left(\V_m b_m\right) &=& \V_m j_R(a_m)\V_m j_R(b_m) = \V_m\left(\V_{pq=m} j_R(a_p) j_R(b_q)\right)\\
&\le&\V_m\left(\V_{pq=m} j_R(a_p b_q)\right)\le\V_m j_R\left(\V_{pq=m} a_p b_q\right)\\
&=&k\left(\V_m\left(\V_{pq=m} a_p b_q\right)\right)=k\left(\V_m a_m \V_m b_m\right)\;.
\end{eqnarray*}
The condition relating $k$ to the involution is equally simple and we omit it. Finally,
it is obvious that for $(a,b)\in R$ we have $k(a)\le k(b)$, since $j_R(a)\le j_R(b)$. But, by definition, $j_R$ is the least nucleus that satisfies this condition, and therefore we conclude that $k=j_R$. Hence, $j_R$ is graded. \qed
\end{proof}

\section{Construction of the Lindenbaum quantales}\label{sec:lind}

\paragraph{The involutive tensor quantale of a frame.}

Let $L$ be a frame, and denote by $I$ the free involutive monoid on one generator $\alpha$, 
whose words are the finite sequences of $\alpha$ and $\alpha^\inv$, and whose unit we shall denote by $\varepsilon$.

For each $w\in I$ we shall denote by $L^{(w)}$ the sup-lattice
$L^{\otimes(\vert w\vert +1)}$, where $\vert w\vert$ is the length of the word $w$ (notice the difference with respect to \P\ref{exm:tensor}):
\[
\begin{array}{rcl}
L^{(\varepsilon)}&=&L\;,\\
L^{(\alpha)}=L^{(\alpha^\inv)}&=&L\otimes L\;,\\
L^{(\alpha\alpha)}=L^{(\alpha\alpha^\inv)}=L^{(\alpha^\inv\alpha)}=L^{(\alpha^\inv\alpha^\inv)}&=&L\otimes L\otimes L\;,\\
&\textrm{etc.}
\end{array}
\]
For each $w,w'\in I$ we define a map
\[\varphi_{w,w'}:L^{(w)}\times L^{(w')}\to L^{(ww')}\]
by
\[\varphi_{w,w'}(x_1\otimes\cdots \otimes x_n,y_1\otimes\cdots\otimes y_m)=x_1\otimes\cdots\otimes(x_n\wedge y_1)\otimes\cdots\otimes y_m\;.\]
It is easy to see that this preserves joins in each variable, and thus it defines a sup-lattice homomorphism
\[\overline\varphi_{w,w'}:L^{(w)}\otimes L^{(w')}\to L^{(ww')}\;.\]
Defining
\[\tensor(L)=\bigoplus_{w\in I} L^{(w)}\]
(not the same as $TL$, \cf\ \P\P \ref{exm:tensor} and \ref{exm:graded}),
and using the distributivity of $\otimes$ over $\bigoplus$,
we obtain the following sup-lattice homomorphism
$\tensor(L)\otimes\tensor(L)\to\tensor(L)$,
\[\xymatrix{\tensor(L)\otimes\tensor(L)\ar[r]^-{\cong}&\bigoplus\limits_{w,w'} L^{(w)}\otimes L^{(w')}\ar[rr]^-{\bigoplus\overline\varphi_{w,w'}}&&
\bigoplus\limits_{w,w'} L^{(ww')}\ar[r]&\tensor(L)}\;,\]
where the rightmost homomorphism is the copairing of the family of sup-lattice embeddings $L^{(ww')}\to\tensor(L)$ which is given by the universal property of the sup-lattice coproduct $\bigoplus L^{(ww')}$. Hence, there is a bilinear multiplication $\tensor(L)\times\tensor(L)\to\tensor(L)$. It is defined on pure tensors by
\[(x_1\otimes\cdots \otimes x_n)(y_1\otimes\cdots\otimes y_m)=x_1\otimes\cdots\otimes(x_n\wedge y_1)\otimes\cdots\otimes y_m\;,\]
and it is straightforward to see that it is associative, hence giving us a quantale multiplication on $\tensor(L)$. This multiplication has a unit, which coincides with $1_L\in L^{(\varepsilon)}=L$, and an involution
\[(-)^\inv:\tensor(L)\to\tensor(L)\]
is obtained from the isomorphisms $L^{(w)}\stackrel\cong\to L^{(w^\inv)}$ that are given by
\[x_1\otimes\cdots\otimes x_n\mapsto x_n\otimes\cdots\otimes x_1\;.\]
Hence, $\tensor(L)$ is a unital involutive quantale, and it is clearly graded over $I$, so that we can define:

\begin{definition}
The \emph{tensor involutive quantale of} $L$ is the graded unital involutive quantale $\tensor(L)$.
\end{definition}

There is an obvious homomorphism of involutive monoids $\overline{(-)}:I\to\tensor(L)$ that sends each word $w$ to
$1_L\otimes\cdots\otimes 1_L\in L^{(w)}$. Hence, in particular, $\overline\alpha=1_L\otimes 1_L\in L^{(\alpha)}$ and $\overline{\alpha^\inv}=1_L\otimes 1_L\in L^{(\alpha^\inv)}$.

The quantale $\tensor(L)$ has the following universal property:

\begin{proposition}
Let $Q$ be a unital involutive quantale such that $b=b^2=b^\inv$ for all $b\in \downsegment e_Q$ (in particular, this implies that $\downsegment e_Q$ is a locale). Let also \[h:L\to\downsegment e_Q\] be a homomorphism of locales, and let $a\in Q$. Then there is exactly one homomorphism of unital involutive quantales
\[\vartheta:\tensor(L)\to Q\]
such that:
\begin{enumerate}
\item\label{univtens1} $\vartheta(x)=h(x)$ for all $x\in L^{(\varepsilon)}=L$;
\item\label{univtens2} $\vartheta(\overline\alpha)=a$.
\end{enumerate}
\end{proposition}

\begin{proof}
By the universal property of the coproduct of sup-lattices, every sup-lattice homomorphism $\vartheta$ with domain $\tensor(L)$ is uniquely determined by its value on the pure tensors of $\tensor(L)$. Furthermore,  every pure tensor \[x_1\otimes\cdots\otimes x_n \in L^{w}\] with $w_i=\alpha,\alpha^\inv$ can be written as a product of $x_i \in L$, $\overline\alpha\in L^{\alpha}$ and $\overline\alpha^\inv\in L^{\alpha^\inv}$, and  thus if $\vartheta$ is a homomorphism of unital involutive quantales that satisfies \ref{univtens1} and \ref{univtens2} then its value is uniquely determined on all the pure tensors. This proves that if $\vartheta$ exists then it is unique. In order to prove existence, assign to each pure tensor
 \[x_1\otimes\cdots\otimes x_n \in L^{w}\]
the value $h(x_1)a_1 h(x_2)\ldots a_{n-1} h(x_n)\in Q$ where $a_i=a$ or $a_i=a^\inv$ according to whether $w_i$ is $\alpha$ or $\alpha^\inv$, respectively. This assignment preserves joins in each variable $x_i$ separately and thus it defines a sup-lattice homomorphism $\vartheta_w:L^{(w)}\to Q$.  The copairing
\[\vartheta=[\vartheta_w]_{w\in I}:\tensor(L)\to Q\]
is easily seen to preserve the quantale multiplication, the unit and the involution, and it satisfies conditions \ref{univtens1} and \ref{univtens2}. \qed
\end{proof}

For the following we denote by $\frames$ the full subcategory of $\uiq$ whose objects are the locales (this is usually called the category of \emph{frames} \cite{Johnstone}).

\begin{corollary}
Let $\wellsppqu$ be the full subcategory of $\uiq$ whose objects are those unital involutive quantales $Q$ such that $b=b^2=b^\inv$ for all $b\le e$ in $Q$. 
Let also $\wellsppqu_*$ be the corresponding category of pointed quantales. 
There is an obvious functor $\wellsppqu_*\to\frames$ that to each quantale $Q$ assigns $\downsegment e$, and such that $h\mapsto h\vert_{\downsegment e}$ for each homomorphism $h$. 
This functor has a left adjoint which to each locale $L$ assigns the pointed quantale $(\tensor(L),\overline\alpha)$.
\end{corollary}

\paragraph{Bimodal frames and pointed quantales.}

We have already mentioned that for an ssq $Q$, $\spp Q = \downsegment e$ is a locale. 
It is also clear that, for $\alpha$ in $Q$, the operators $\lozenge x=\spp(\alpha x)$ and $\blacklozenge x=\spp(\alpha^\inv x )$ preserve arbitrary joins, hence they are sup-lattice endomorphisms of $\spp Q$.

\begin{definition}\label{def:conjug}
 We say that two sup-lattice endomorphisms of $L$, $\lozenge$ and $\blacklozenge$, are \emph{conjugate modalities} if for all $x,y \in L$ we have
\begin{eqnarray*}
\lozenge x\wedge y & \le & \lozenge (x\wedge\blacklozenge y)\;,\\
\blacklozenge x\wedge y & \le & \blacklozenge (x\wedge\lozenge y)\;.
\end{eqnarray*}
A \emph{bimodal frame} $(L,\lozenge,\blacklozenge)$ is a frame $L$ equipped with two conjugate modalities $\lozenge$ and $\blacklozenge$.\\
\end{definition}

\begin{proposition}
Let $Q$ be an ssq, and $\alpha\in Q$. Then \[(\spp Q,\spp(\alpha \, -),\spp(\alpha^\inv \, -))\] is a bimodal frame.
\end{proposition}

\begin{proof}
We only have to check that the two endomorphisms are conjugate modalities. From the fact that $\spp$ is a support (\ref{spp:2}) we have:
\[\spp(\alpha x)y\le\alpha x \alpha^\inv y\;,\] and thus 
\[\begin{array}{rcl}&\lozenge x\wedge y = \spp(\alpha x)y= \spp(\spp(\alpha x)y)
\le\spp(\alpha x \alpha^\inv y)
=\spp(\alpha x \spp(\alpha^\inv y))=\lozenge(x\wedge\blacklozenge y)\;,
\end{array}\] 
using stability and the fact that $x,y \in \spp Q$. The other conjugacy condition is obtained by interchanging $\alpha$ and $\alpha^\inv$. \qed
\end{proof}

From now on we shall refer to any sup-lattice endomorphism
$\spp:Q\to Q$
on a unital involutive quantale  $Q$
such that $\spp a\le e$ for all $a\in Q$ as a \emph{pre-support} of $Q$.

Given a bimodal frame $(L,\lozenge,\blacklozenge)$, a pre-support can be easily defined on pure tensors of $\tensor(L)$ (and then extended to joins of these in the obvious way): if $x=x_0\otimes\cdots\otimes x_n$ is in degree
$w=w_1\ldots w_n$ (with  $w_i\in\{\alpha,\alpha^\inv\}$) then
\[\spp x=x_0\wedge \langle w_1\rangle(x_1\wedge\langle w_2\rangle(\ldots))\;,\]
where $\langle w_i\rangle$ is $\lozenge$ or $\blacklozenge$ according to whether $w_i=\alpha$ or $w_i=\alpha^\inv$, respectively.
Recursively, we have:

\begin{definition}\label{suptensor}
Let $x=x_0\otimes\cdots\otimes x_n\in L^{(w_1\ldots w_n)}$. Then,
\begin{itemize}
\item $\spp x=x$, if $n=0$;
\item $\spp x=x_0\wedge\langle w_1\rangle(\spp x')$, if $n\ge 1$, where $x'=x_1\otimes\cdots \otimes x_n\in L^{(w_2\ldots w_n)}$.
\end{itemize}
\end{definition}

\begin{lemma}\label{lemmaA}
The following properties hold for all $a,b,c\in\tensor(L)$:
\begin{enumerate}
\item $\spp e=e$ and $\spp a \le e$ (condition \ref{spp:1} in the definition of support);
\item\label{outralinea} $\spp(\spp a b)=\spp a\spp b$ (in particular, $\spp\spp a=\spp a$);
\item $\spp(ab)=\spp(a\spp b)$ (we say the pre-support $\spp$ is \emph{stable});
\item If $\lozenge$ and $\blacklozenge$ are conjugate we have
\begin{enumerate}
\item\label{ala} $\spp a \le\spp(aa^\inv)$;
\item\label{alb} $\spp (\spp a b) \le\spp(aa^\inv b)$;
\item\label{alc} $\spp(c\spp ab)\le\spp(caa^\inv b)$.
\end{enumerate}
\end{enumerate}
\end{lemma}
\begin{proof}
\begin{enumerate}
	\item $\spp \tensor(L) = L$ and $e$ is $1_L$.
	\item  It suffices to prove this for all the pure tensors $b\in L^{(v)}$, where $v$ is an arbitrary degree. Let
	\[
	b=y_0\otimes\cdots\otimes y_p\in L^{(v_1\ldots v_p)}\;.
\]
We have $\spp a\in L^{(\varepsilon)}$, and thus from \P\ref{suptensor} we obtain 
\begin{eqnarray*}
\spp(\spp a b)&=& \spp((\spp a\wedge y_0)\otimes\cdots\otimes y_p)\\
	&=& (\spp a \wedge y_0)\wedge \op{v_1} (\spp(y_1\otimes\cdots\otimes y_p))\\
	&=& \spp a \wedge (y_0\wedge \op{v_1}(\spp(y_1\otimes\cdots\otimes y_p)))\\
	&=& \spp a\spp b\;.
\end{eqnarray*}
	\item It suffices to prove this for all the pure tensors $a\in L^{(w)}$, where $w$ is an arbitrary degree. Let then
\[
	a=x_0\otimes\cdots\otimes x_n\in L^{(w_1\ldots w_n)}\;.
\]
The proof is done by induction on $n$. For the base case assume that $n=0$; that is, we have $a=\spp a=x_0\in L^{(\varepsilon)}$, and from \ref{outralinea} we obtain 
\[
\spp(a b)= \spp(\spp a b) =\spp a\spp b=a\spp b\;,
\]
whence $\spp(ab)=\spp\spp(ab)=\spp(a\spp b)$.

Now for the induction step let $n\ge 1$ and let
\[
a'=x_0\otimes\cdots\otimes x_{n-1}\in L^{(w_1\ldots w_{n-1})}\;.
\]
We have
\[
a=a'\overline{w_n}x_{n}
\]
and thus
\[\begin{array}{rcll}
	\spp(a  b)&=&\spp(a'\overline{w_n}x_{n}b)\\
	&=&\spp(a'   \spp(\overline{w_n}x_{n}b)) & \textrm{(Induction hyp.)}\\
	&=&\spp(a'  \op{w_{n}}( \spp (x_n b))) & \textrm{(Def.\ of }\spp\textrm{)}\\
	&=&\spp(a'  \op{w_{n}}(x_n \spp  b)) & \textrm{(By \ref{outralinea})}\\
	&=&\spp(a'  \spp(\overline{w_n}x_n \spp b))) & \textrm{(Def.\ of }\spp\textrm{)}\\
	&=&\spp(a'  \overline{w_n}x_n \spp b) & \textrm{(Induction hyp.)}\\
	&=& \spp(a \spp b) 	\;.
	\end{array}\]
         \item First we remark that (\ref{ala}) is an instance of (\ref{alb}) (which in turn is an instance of (\ref{alc})).
         Moreover, (\ref{alb}) implies (\ref{alc}) due to stability: if we assume
         (\ref{alb}) then
         \[\spp(c\spp ab)=\spp(c\spp(\spp ab))\le\spp(c\spp(aa^\inv b))=\spp(caa^\inv b)\;.\]
         It now suffices to prove (\ref{alb}). We shall prove $\spp(\spp a b)\le\spp(aa^\inv b)$ for the particular case
         where $a$ is a pure tensor
         \[a=x_0\otimes\cdots\otimes x_n\]
         in degree $w=w_1\ldots w_n$ ($n\ge 1$),
         which implies the general case. The proof is by induction on $n$.
	 
	 Base: for $n=0$, we have $a=\spp a=aa^\inv$, and thus
	 $\spp(\spp ab)=\spp(aa^\inv b)$.
	 
	 Step: for $n\ge 1$, let $r$ be $x_1\otimes\cdots\otimes x_n \in L^{(w_2\ldots w_n)}$; hence, we have $a=x_0\overline{w_1} r$, and thus using stability we obtain
\[\spp(aa^\inv b)=\spp(x_0\overline{w_1} rr^\inv\overline{w_1^\inv} x_0 b)
                     =\spp(x_0\overline{w_1} rr^\inv\spp(\overline{w_1^\inv} x_0 \spp b))\;.\]
By the definition of $\spp$ (and making the symbol $\wedge$ explicit for the multiplication in $\spp\tensor(L)$) this equals
\[x_0\wedge \op{w_1}\spp(rr^\inv \op{w_1^\inv}(x_0 \wedge\spp b))\;,\]
which, by the induction hypothesis, is greater or equal to
\[x_0\wedge \op{w_1}\spp(\spp r\wedge \op{w_1^\inv}(x_0 \wedge\spp b))\;,\]
which in turn equals
\[x_0\wedge \op{w_1}(\spp r\wedge \op{w_1^\inv}(x_0 \wedge\spp b))\]
(because the argument of the outermost occurrence of $\spp$ was in $\spp\tensor(L)$).
Finally, by conjugacy of the operators $\op{w_1}$ and $\op{w_1^\inv}$ the latter is greater or equal to
\begin{eqnarray*}
x_0\wedge \op{w_1}(\spp r)\wedge x_0 \wedge\spp b&=&
x_0\wedge \op{w_1}(\spp r) \wedge\spp b\\
&=&\spp(x_0\overline{w_1} r)\wedge \spp b=\spp a\wedge \spp b=\spp(\spp a b)\;. \qed
\end{eqnarray*}
\end{enumerate}

\end{proof}

\paragraph{The supported quantale of a bimodal frame.}

So far we have obtained, from an arbitrary bimodal frame, a quantale with a stable pre-support. In order to obtain an actual supported quantale we shall impose the missing properties, namely $\spp a \le aa^\inv$ (\ref{spp:2}) and $a \le \spp aa$ (\ref{spp:3}), by taking quotients of $\tensor(L)$.

\begin{definition} \label{tensorcomsup}
Let $j_\spp$ be the least nucleus $j$ on $\tensor(L)$ such that
\[j(a)=j(\spp a a)\;.\] We define
$\sppT(L,\lozenge,\blacklozenge)$ to be $\tensor(L)_{j_\spp}$.
We also write $\sppT(L)$ if $\lozenge$ and $\blacklozenge$ are clear from the context.
\end{definition}

From \P\ref{prop:gradedrelation} it is clear that $j_\spp$ is a graded nucleus, and it is the identity on $L^{(\varepsilon)}=L$ because for all $a \in L$ we have  $\spp a a = a$. Hence, we have concluded, just from the graded structure of $\tensor(L)$, that the injection of generators of $L$ into $\sppT(L,\lozenge,\blacklozenge)$ is 1--1.

Proving a similar fact for the other axiom, $\spp a\le aa^\inv$, is less easy, and we shall address this now.
Let $(L,\lozenge,\blacklozenge)$ be a bimodal frame and let \[R=\{(\spp a,aa^\inv) \mid a \in \sppT(L,\lozenge,\blacklozenge)\}\;.\]

\begin{definition}  
$\K(L,\lozenge,\blacklozenge)$ is $\sppT(L,\lozenge,\blacklozenge)_{j_{R}}$. As in \P\ref{tensorcomsup} we may write $\K(L)$. We shall denote the selected point $j_R(j_\spp(\overline\alpha))\in \K(L)$ by $\boldsymbol{\alpha}$.
\end{definition}

\begin{lemma}\label{lemmaB}
Recall the definition (\ref{Rbarra}) of $\overline R$.
If $(y,z)\in\overline R$ then $\spp y\le\spp z$ in $\sppT(L,\lozenge,\blacklozenge)$.
\end{lemma} 

\begin{proof}
We have $\spp y\le\spp z$ for all $(y,z)\in\overline R$ if and only if the following two conditions hold for all $(y,z)\in\overline R$ and all $a,b\in\sppT(L)$:
\begin{eqnarray}
\spp(ayb)&\le&\spp(azb) \label{prop1}\\
\spp(ay^\inv b)&\le&\spp(az^\inv b)\;.\label{prop2}
\end{eqnarray}
In order to prove these two conditions we shall show that they hold for all $(y,z)\in R$ and that they are preserved by the recursive rules of construction of $\overline R$.

Let $(y,z)\in R$; that is, let $y$ and $z$ be of the form $y=\spp t$ and $z=tt^\inv$. We have
\[\spp(ayb)=\spp(a\spp t b)\le\spp(att^\inv b)=\spp(azb)\;,\]
from \P\ref{lemmaA}. We similarly have
$\spp(ay^\inv b)\le\spp(az^\inv b)$ because $y$ and $z$ are self-adjoint.

Now assume that (\ref{prop1})--(\ref{prop2}) hold for some pair $(y,z)\in\sppT(L)\times\sppT(L)$. We shall prove that they equally hold for the following pairs: (i) $(\spp y,\spp z)$; (ii) $(y^\inv,z^\inv)$; (iii) $(qy,qz)$; and (iv) $(yq,zq)$, for all $q\in\sppT(L)$.

(i) Since $(y,z)$ satisfies (\ref{prop1})--(\ref{prop2}) we have $\spp y\le\spp z$ (make $a=b=1_L$), and thus
$\spp(a\spp y b)\le\spp(a\spp z b)$ for all $a,b\in\sppT(L)$, proving (\ref{prop1}) for the pair $(\spp y,\spp t)$. Since $\spp y$ and $\spp z$ are self-adjoint, we also conclude (\ref{prop2}) for the pair $(\spp y,\spp t)$.

(ii), (iii) and (iv) are obvious. \qed
\end{proof}

\begin{theorem}\label{theorem}
The unit of the adjunction $L\to \K(L)$ is 1--1.
\end{theorem}

\begin{proof}
Equivalently, we want to  prove that $\spp\sppT(L)\subset Q_{j_{R}}$, or, in other words,
that for all $(y,z)\in\overline R$ and $x\in L$ we have
\begin{equation}
z\le x\Rightarrow y\le x\;.\label{Pcondition}
\end{equation}
Let $P\subset\overline R$ be the subset of $\overline R$ consisting of all those $(y,z)$ such that (\ref{Pcondition}) holds for all $x\in L$. We shall prove that $R\subset P$ and that $P$ is closed under the recursive formation rules of $\overline R$, hence showing that $P=\overline R$ because $\overline R$ is the least subset of $\sppT(L)\times\sppT(L)$ satisfying these conditions.

Let $(\spp y,yy^\inv)\in R$. If $x\in L$ and $yy^\inv\le x$ we conclude that $y=\spp y=yy^\inv\in L$ due to the grading of $\sppT(L)$ over the free involutive monoid $I$, and thus $\spp y\le x$, showing that $R\subset P$.

From now on let $(y,z)$ be a fixed but arbitrary element of $P$. By \P\ref{lemmaB} we conclude $\spp y\le\spp z$. Hence, (\ref{Pcondition}) holds for the pair $(\spp y,\spp z)$, and thus $(\spp y,\spp z)\in P$.
Now let $q\in\sppT(L)$ and assume that $qz\le x$ for some $x\in L$. Then $z\le 1_L$ (again due to the grading over $I$), and thus $y\le 1_L$ because $(y,z)\in P$ (make $x=1_L$). Hence, again using the previous lemma we obtain
\[qy=q\spp y\le q\spp z= qz\le x\;,\]
showing that $(qy,qz)\in P$. In a similar way we conclude that $(yq,zq)\in P$. Finally, $x\le 1_L$ implies that $x$ is self-adjoint, and thus the conditions $y^\inv\le x$ and $z^\inv\le x$ are equivalent to $y\le x$ and $z\le x$, respectively, showing that $(y^\inv,z^\inv)\in P$. \qed
\end{proof}

\paragraph{T, K4, S4, S5.}

Now let us extend our results to the systems of modal logic T, K4, S4, and S5. As was explained in \S\ref{logica}, we shall need to impose additional conditions on the selected element $\boldsymbol{\alpha}\in\K$, such as reflexivity ($\boldsymbol{\alpha}\ge 1_L$), transitivity ($\boldsymbol{\alpha}^2\le\boldsymbol{\alpha}$), etc. In order to obtain again coreflections we shall also have to impose additional axioms on the modalities $\lozenge$ and $\blacklozenge$ of $L$. As we mentioned in \S\ref{introduction}, we shall see that for T, K4 and S4 these coincide with the well known axioms for the corresponding systems of modal logic under consideration (in other words, the same axioms still produce complete axiomatizations for the new semantics), whereas for S5 a new axiomatization is defined just by imposing that $\lozenge$ should coincide with $\blacklozenge$.

The proof techniques are very similar to those employed in the previous sections for the system K. In fact we could have already presented the theory for K in such a generality so as to be able to directly reuse the results now, but this would have obscured the main ideas, so for the sake of clarity we shall single out the general results only now.

\begin{lemma}\label{lemmaBnew}
Let $\rho\subset\sppT(L)\times\sppT(L)$ be any binary relation on $\sppT(L)$, and let $\overline\rho$ be the closure of $\rho$ under the rules
\begin{eqnarray*}
(y,z)\in\overline\rho&\Rightarrow& (\spp y,\spp z)\in\overline\rho\\
(y,z)\in\overline\rho&\Rightarrow&(ayb,azb)\in\overline\rho\textrm{ for all }a,b\in\sppT(L)\\
(y,z)\in\overline\rho&\Rightarrow& (y^\inv,z^\inv)\in\overline\rho\;.
\end{eqnarray*}
Assume that for all $(y,z)\in\rho$ and all $a,b\in\sppT(L)$ we have
\begin{eqnarray*}
\spp(ayb)&\le&\spp(azb)\\
\spp(ay^\inv b)&\le&\spp(az^\inv b)\;.
\end{eqnarray*}
Then for all $(y,z)\in\overline\rho$ we have $\spp y\le\spp z$.
\end{lemma}

\begin{proof}
This follows from a simple adaptation of the proof of \P\ref{lemmaB}. \qed
\end{proof}

\begin{lemma}\label{lemmaC}
Let $P\subset\sppT(L)\times\sppT(L)$ be the set of all those $(y,z)$ satisfying the following two conditions:
\begin{enumerate}
\item $\spp y\le\spp z$;
\item $z\le x\Rightarrow y\le x$ for all $x\in L$.
\end{enumerate}
Then $P$ is closed under the rules
\begin{eqnarray*}
(y,z)\in P&\Rightarrow& (\spp y,\spp z)\in P\\
(y,z)\in P&\Rightarrow&(ayb,azb)\in P\textrm{ for all }a,b\in\sppT(L)\\
(y,z)\in P&\Rightarrow& (y^\inv,z^\inv)\in P\;.
\end{eqnarray*}
\end{lemma}

\begin{proof}
The proof of this lemma is contained in the proof of \P\ref{theorem}, where $P$ was defined to be a subset of $\overline R$, but in fact the only property of $\overline R$ used in order to prove the closure properties of $P$ was the fact that for all $(y,z)\in\overline R$ we have $\spp y\le\spp z$. The other key ingredient is the fact that $ab\in L$ implies $a,b\in L$ for all $a,b\in\sppT(L)$, due to the grading of $\sppT(L)$ over $I$. \qed
\end{proof}

\begin{theorem}
Let $\rho\subset\sppT(L)\times\sppT(L)$ be a binary relation such that for all $(y,z)\in\rho$ and all $a,b\in\sppT(L)$ we have
\begin{eqnarray}
\spp(ayb)&\le&\spp(azb) \label{relone}\\
\spp(ay^\inv b)&\le&\spp(az^\inv b)\label{reltwo}\;,
\end{eqnarray}
and let $Q$ be the (supported) quotient of $\sppT(L)$ generated by $\rho$. Then the injection of generators of $L$ onto $\spp Q$,
\[\eta:L\to\spp\sppT(L)\to\spp Q\;,\]
is an isomorphism.
\end{theorem}

\begin{proof}
This is an immediate consequence of the previous two lemmas, by a reasoning analogous to that of \P\ref{theorem}. \qed
\end{proof}

Similarly to what we have done in \S\ref{logica} for the Lindenbaum quantales $\LindT$, $\LindKfour$, etc.\ (see \P\ref{def:Lindquantales}),
we define $\TT(L)$, $\TKq(L)$,$\TSq(L)$ and $\TSc(L)$ to be quotients of $\K(L)$ by analogous defining relations:
\begin{description}
\item[$\TT(L)$:] $e\le\boldsymbol{\alpha}$
\item[$\TKq(L)$:] $\boldsymbol{\alpha}\boldsymbol{\alpha}\le\boldsymbol{\alpha}$
\item[$\TSq(L)$:] $e\le\boldsymbol{\alpha} \ge\boldsymbol{\alpha}\boldsymbol{\alpha}$
\item[$\TSc(L)$:] $e\le\boldsymbol{\alpha}=\boldsymbol{\alpha}^\inv\ge\boldsymbol{\alpha}\boldsymbol{\alpha}$
\end{description}

\begin{corollary}\label{cor:TT}
Let $L$ be a bimodal frame such that for all $x\in L$ the conditions $x\le\lozenge x$ and $x\le\blacklozenge x$ hold. Then the injection of generators $L\to\spp\TT(L)$ is an isomorphism.
\end{corollary}

\begin{proof}
The quantale $\TT(L)$ is the quotient of $\sppT(L)$ generated by the conditions
\begin{eqnarray*}	
\spp  a  &\le& aa^\inv, \textrm{ for all } a \in \sppT(L)\,,\\
\boldsymbol{\alpha}&\ge& e\,.
\end{eqnarray*}
Define $\RT\subset\sppT(L)\times\sppT(L)$ as follows:
\[\RT=\{(\spp y,yy^\inv)\st y\in\sppT(L)\})\cup\{(1_L,\boldsymbol{\alpha})\}\;.\]
All we have to do is, by the previous theorem, prove that for all $(y,z)\in\RT$ the conditions (\ref{relone})--(\ref{reltwo})
are satisfied for all $a,b\in\sppT(L)$. This has already been done for the pairs of the form
$(\spp y, yy^\inv)$ in \P\ref{lemmaB}, so we only have to concern ourselves with the pair $(1_L,\boldsymbol{\alpha})$.
For $a,b \in Q$ we have
\[
\spp(ab)=\spp(a \spp(b))\le \spp(a \lozenge (\spp(b))) = \spp(a \spp(\boldsymbol{\alpha} b)) = \spp(a \boldsymbol{\alpha} b)\;,
\]
where we have used stability of $\spp$ twice and the inequality follows from $x\le\lozenge x$ and monotonicity of $\spp$; this proves (\ref{relone}).
Then (\ref{reltwo}) is proved in a similar way using the inequality $x\le\blacklozenge x$. \qed
\end{proof}

\begin{corollary}\label{cor:K4}
Let $L$ be a bimodal frame such that for all $x\in L$ the conditions
$\lozenge \lozenge(x) \le \lozenge(x)$ and $\blacklozenge\blacklozenge(x) \le \blacklozenge(x)$ hold. 
Then the injection of generators $L\to\TKq(L)$ is an isomorphism.
\end{corollary}

\begin{proof}
It remains to prove that the pair
$(y,z)=(\boldsymbol{\alpha}^2,\boldsymbol{\alpha})$ satisfies (\ref{relone})--(\ref{reltwo}).
The first condition is proved as follows:
\[
\spp(a \boldsymbol{\alpha} b)=\spp(a\spp(\boldsymbol{\alpha} b))=\spp(a \lozenge (\spp(b))) \ge \spp(a \lozenge^2 (\spp(b))) = \spp(a \boldsymbol{\alpha}^2 b)\,.\]
(\ref{reltwo}) is proved in the same way once we replace $\boldsymbol{\alpha}$ by $\boldsymbol{\alpha}^\inv$ and $\lozenge$ by $\blacklozenge$ in the previous argument.
\qed
\end{proof}

\begin{corollary}\label{cor:S4}
Let $L$ be a bimodal frame such that for all $x\in L$ the conditions $x\le\lozenge x$, $x\le\blacklozenge x$, $\lozenge \lozenge(x) \le \lozenge(x)$ and $\blacklozenge\blacklozenge(x) \le \blacklozenge(x)$ hold. Then the injection of generators $L\to\TSq(L)$ is an isomorphism.
\end{corollary}

\begin{proof}
Immediate from the previous two corollaries. \qed
\end{proof}

\begin{corollary}\label{cor:S5}
Let $L$ be a bimodal frame such that for all $x\in L$ the conditions of the previous corollary and $\lozenge x = \blacklozenge x$ hold. Then the injection of generators $L\to\TSc(L)$ is an isomorphism.
\end{corollary}

\begin{proof}
It remains to show that the pairs $(\boldsymbol{\alpha}, \boldsymbol{\alpha}^\inv)$ and $(\boldsymbol{\alpha}^\inv, \boldsymbol{\alpha})$ satisfy (\ref{relone}) and (\ref{reltwo}), which  is done as follows:
\[
\spp(a \boldsymbol{\alpha} b)=\spp(a\spp(\boldsymbol{\alpha} b))=\spp(a \lozenge (\spp(b))) = \spp(a \blacklozenge (\spp(b))) = \spp(a \boldsymbol{\alpha}^\inv b)\;. \qed
\]
\end{proof}

\paragraph{The Lindenbaum quantales.}

Let $\BK$ be the Lindenbaum algebra for system K, as in \S\ref{logica}, and let $\BT$, $\BKq$, $\BSq$, and $\BSc$ be the Lindenbaum algebras for systems T, K4, S4, and S5, respectively. These are \emph{modal lattices} in the following sense:

\begin{definition}
By a \emph{modal lattice} is meant a bounded distributive lattice $L$ equipped with an endomap $\lozenge$ that preserves finite joins. A \emph{bimodal lattice} is a modal lattice equipped with another endomap $\blacklozenge$ that preserves finite joins and  in addition satisfies conjugacy relations similar to those of bimodal frames:
\begin{eqnarray*}
\lozenge x\land y & \le & \lozenge (x \land \blacklozenge y)\;,\\
\blacklozenge x\land y & \le & \blacklozenge (x \land \lozenge y)\;.\\
\end{eqnarray*}
The \emph{category $\modlatcat$ of modal lattices} has the modal lattices as objects and the homomorphisms of bounded lattices that preserve $\lozenge$ as morphisms. The \emph{category $\bimodlatcat$ of bimodal lattices} is defined analogously, with objects being the bimodal lattices and the morphisms being the homomorphisms of modal lattices that also preserve $\blacklozenge$.

We shall also refer to a modal lattice as a
\begin{itemize}
\item {\it T-modal lattice} if $x\le\lozenge x$ for all $x$;
\item {\it K4-modal lattice} if $\lozenge\lozenge x\le\lozenge x$ for all $x$;
\item {\it S4-modal lattice} if it is both a T-modal lattice and a K4-modal lattice;
\item {\it S5-modal lattice} if it is an S4-modal lattice and $\lozenge x\land y\le\lozenge(x\land\lozenge y)$ for all $x$ and $y$ (\ie, $\lozenge$ is conjugate to itself).
\end{itemize}
The categories $\Tmodlatcat$, $\Kqmodlatcat$, $\Sqmodlatcat$, and $\Scmodlatcat$ are, respectively, the full subcategories of $\modlatcat$ whose objects are the T-modal lattices, the K4-modal lattices, the S4-modal lattices, and the S5-modal lattices.

For bimodal lattices we adopt a similar terminology: a bimodal lattice is referred to as a
\begin{itemize}
\item {\it T-bimodal lattice} if $x\le\lozenge x$ and $x\le\blacklozenge x$ for all $x$;
\item {\it K4-bimodal lattice} if $\lozenge\lozenge x\le\lozenge x$ and $\blacklozenge\blacklozenge x\le\blacklozenge x$ for all $x$;
\item {\it S4-bimodal lattice} if it is both a T-bimodal lattice and a K4-bimodal lattice;
\item {\it S5-bimodal lattice} if it is an S4-bimodal lattice and $\lozenge x=\blacklozenge x$ for all $x$.
\end{itemize}
The categories $\Tbimodlatcat$, $\Kqbimodlatcat$, $\Sqbimodlatcat$, and $\Scbimodlatcat$ are, respectively, the full subcategories of $\bimodlatcat$ whose objects are the T-bimodal lattices, the K4-bimodal lattices, the S4-bimodal lattices, and the S5-bimodal lattices.
\end{definition}

By standard universal algebra the forgetful functor $\bimodlatcat\to\modlatcat$ has a left adjoint which assigns to $\BK$ its ``enveloping'' bimodal lattice $\BK'$. Similarly, there are left adjoints
\begin{eqnarray*}
\Tmodlatcat&\to&\Tbimodlatcat\\
\Kqmodlatcat&\to&\Kqbimodlatcat\\
\Sqmodlatcat&\to&\Sqbimodlatcat\;,
\end{eqnarray*}
and we write $\BT'$, $\BKq'$, and $\BSq'$ for the respective  images of
$\BT$, $\BKq$, and $\BSq$ under the left adjoints.
For S5 the situation is simpler because the categories $\Scmodlatcat$ and $\Scbimodlatcat$ are obviously isomorphic, since any S5-modal lattice becomes an S5-bimodal lattice just by defining $\blacklozenge$ to coincide with $\lozenge$. We shall also write $\BSc'$ for $\BSc$ thus regarded as a bimodal lattice.

Since $\lozenge$ and $\blacklozenge$ preserve finite joins they can be extended canonically to sup-lattice endomorphisms of the ideal completion
$\Idl(\BK')$, which is a frame because $\BK'$ is a distributive lattice. The conjugation relations are easily seen to be inherited from those of $\BK'$, and thus $\Idl(\BK')$ is a bimodal frame. Similar remarks apply to the other Lindenbaum algebras, and in addition $\Idl(\BT')$ satisfies the axioms of a T-bimodal lattice, $\Idl(\BKq')$ satisfies the axioms of a K4-bimodal lattice, etc. (Hence, in particular, the propositions in \P\P \ref{cor:TT}--\ref{cor:S5} can be applied to
$\Idl(\BT')$, $\Idl(\BKq')$, $\Idl(\BSq')$, and $\Idl(\BSc')$, respectively.)

Summarizing, we have described a way of constructing functors from modal lattices to bimodal frames which are left adjoint to the obvious forgetful functors. Composing these functors with the left adjoints from bimodal frames to pointed ssqs we obtain from $\BK$, $\BT$, $\BKq$, $\BSq$, and $\BSc$ pointed
ssqs $\K(\Idl(\BK'))$, $\TT(\Idl(\BT'))$, $\TKq(\Idl(\BKq'))$, $\TSq(\Idl(\BSq'))$, and $\TSc(\Idl(\BSc'))$, respectively. For each system $S\in\{\textrm{K},\textrm{T},\textrm{K4},\textrm{S4},\textrm{S5}\}$ the map obtained by composing the following arrows (the leftmost one is just the natural quotient),
\[\BK\to B_S\to B'_S\to \Idl(B'_S)\to \tensor_S(B'_S)\;,\]
has the same universal property as the injection of generators
$\BK\to\Lind_S$.
Hence, the quantale $\tensor_S(B'_S)$ is a particular construction of the Lindenbaum quantale $\Lind_S$:

\begin{theorem}
\begin{eqnarray*}
\LindK&\cong&\K(\Idl(\BK'))\\
\LindT&\cong&\TT(\Idl(\BT'))\\
\LindKfour&\cong&\TKq(\Idl(\BKq'))\\
\LindSfour&\cong&\TSq(\Idl(\BSq'))\\
\LindSfive&\cong&\TSc(\Idl(\BSc'))
\end{eqnarray*}
\end{theorem}

Since, as we have remarked above, S5-modal lattices and S5-bimodal lattices are ``the same'', our results immediately tell us that the unit of the adjunction between S5-modal lattices and pointed ssqs,
\[\BSc\to\LindSfive\;,\]
is a monomorphism. In logical terms this means that a complete axiomatization for the system S5 (with the advantage of making no use of negation or the modal necessity operator) can be as follows:

\begin{theorem}
S5 is complete for the following axiom schemata:
\begin{eqnarray*}
\varphi&\rightarrow&\lozenge\varphi\\
\lozenge\lozenge\varphi&\rightarrow&\lozenge\varphi\\
\lozenge\varphi\land\psi&\rightarrow&\lozenge(\varphi\land\lozenge\psi)\;.
\end{eqnarray*}
\end{theorem}

We have not verified whether the remaining canonical mappings
\begin{eqnarray*}
\BK&\to&\BK'\\
\BT&\to&\BT'\\
\BKq&\to&\BKq'\\
\BSq&\to&\BSq'
\end{eqnarray*}
are monomorphisms (although we believe they are). This means that we have not verified completeness for the classical axiomatizations of K, T, K4, and S4. Of course, by this we mean we have not verified this \emph{in an arbitrary topos}, for otherwise we know, from the classical completeness theorems of propositional normal modal logic, that the axiomatizations are complete: using Zorn's Lemma we can find a Kripke structure $(W,R)$ that gives us a monomorphism
\[\BK\to \BK'\to\K(\Idl(\BK'))\cong\LindK\to\pwset{W\times W}\]
implying that $\BK\to \BK'$ is 1--1, and  the same applies to T, K4, and S4.

~\\
{\sc
Centro de L\'ogica e Computa\c{c}\~ao\\
Instituto Superior T\'{e}cnico\\
Universidade T\'{e}cnica de Lisboa\\
Av.\ Rovisco Pais 1, 1049-001 Lisboa, Portugal}\\
{\it E-mail:} {\sf sergiortm@gmail.com}\\
~\\
{\sc
Centro de An\'alise Matem\'atica, Geometria e Sistemas Din\^amicos
Departamento de Matem\'{a}tica do Instituto Superior T\'{e}cnico\\
Universidade T\'{e}cnica de Lisboa\\
Av.\ Rovisco Pais 1, 1049-001 Lisboa, Portugal}\\
{\it E-mail:} {\sf pmr@math.ist.utl.pt}

\end{document}